\documentclass[a4paper]{article}

\usepackage[english]{babel}
\usepackage[utf8]{inputenc}
\usepackage{amsmath}
\usepackage{amssymb}
\usepackage{amsthm}
\usepackage{graphicx}
\usepackage[colorinlistoftodos]{todonotes}

\theoremstyle{plain}
\newtheorem{defi}{Definition}[section]
\newtheorem{lemma}[defi]{Lemma}
\newtheorem{coro}[defi]{Corollary}
\newtheorem{prop}[defi]{Proposition}
\newtheorem{theo}[defi]{Theorem}
\theoremstyle{definition}
\newtheorem{example}[defi]{Example}
\theoremstyle{remark}
\newtheorem{remark}[defi]{Remark}
\numberwithin{equation}{section}

\def\nn{\mathbb{N}}
\def\zz{\mathbb{Z}}
\def\rr{\mathbb{R}}
\def\cc{\mathbb{C}}

\def\M{\mathfrak{M}}
\def\L{\mathfrak L}
\def\h{\mathfrak h}
\def\H{\mathcal H}
\def\D{\mathcal D}
\def\R{\mathcal R}
\def\pre{\noindent{\bf Proof:}}
\def\fin{$\Box$}
\def\tr{\mathrm{Tr}}
\def\id{\mathrm{Id}}
\def\e{\mathrm{e}}
\def\i{\mathrm{i}}
\def\d{\mathrm{d}}
\makeatletter
\newcommand{\oset}[2]{%
  {\mathop{#2}\limits^{\vbox to -.5\ex@{\kern-\tw@\ex@
   \hbox{\scriptsize #1}\vss}}}}
\makeatother
\def\submod{\oset{{\tiny d}}{-}}
\def\addmod{\oset{{\tiny d}}{+}}

\def\ind{1\hspace{-0.27em}\mathrm{l}}
\def\pp{\mathbb{P}}
\def\ee{\mathbb{E}}
\def\rhoinv{\rho^{\mathrm{inv}}}

\def\eps{\varepsilon}

\def\ket#1{|#1\rangle}
\def\braket#1#2{\langle #1  , #2\rangle}
\def\ketbra#1#2{|#1\rangle\langle#2|}

\title{Homogeneous open quantum random Walks\\on a lattice
}
\author{Raffaella Carbone\\Dipartimento di Matematica dell'Universit\'a di Pavia\\ via Ferrata, 1, 27100 Pavia, Italy\\\texttt{raffaella.carbone@unipv.it}\\ \\Yan Pautrat\\Laboratoire de Math\'ematiques\\
Universit\'e Paris-Sud\\
91405 Orsay Cedex, France\\\texttt{yan.pautrat@math.u-psud.fr}}
\begin{document}
\bibliographystyle{abbrv}
\maketitle

{\noindent\bf Abstract.} We study open quantum random walks (OQRW) for which the underlying graph is a lattice, and the generators of the walk are translation-invariant. Using the results obtained recently in \cite{CP1}, we study the quantum trajectory associated with the OQRW, which is described by a position process and a state process. We obtain a central limit theorem and a large deviation principle for the position process. We study in detail the case of homogeneous OQRWs on a lattice, with internal space $\h=\cc^2$.

\section{Introduction}

Open quantum random walks (OQRW in short) were defined by Attal \textit{et al.} in~\cite{APSS}. 
They seem to be a good quantum analogue of Markov chains, and, as such, are a very promising tool to model many physical problems (see \cite{CP1} and the references therein).
In the paper \cite{CP1}, we described the notions of irreducibility and aperiodicity for OQRWs, and derived, in particular, convergence properties for irreducible, or irreducible and aperiodic, OQRWs. In the same way as for classical Markov chains, those convergence results assumed the existence of an invariant state.

In the present paper we focus on translation-invariant OQRWs on a lattice, which attracted special attention in many recent papers (see, for instance, \cite{AGS,BBT,KY,LardSouza,SP2,XY}).
As we have shown in \cite{CP1}, these OQRWs do not have an invariant state so that most of the convergence results from \cite{CP1} are useless. We will show, however, that there exists an auxiliary map which allows to characterize many properties of the homogeneous open quantum random walk. With the help of the Perron-Frobenius theorem described in \cite{CP1}, we can obtain results about the \textit{quantum trajectory} associated with the process. More precisely, if the quantum trajectory is given by the position and state couple $(X_p,\rho_p)$, we obtain a central limit theorem and a large deviation principle for the position process $(X_p)_p$.

The immediate physical application of our results is quantum measurements, or more precisely, repeated indirect measurements. In that framework, a given system $\mathcal S$ interacts sequentially with external systems $\mathcal R_p$, $p=1,2,\ldots$ representing measuring devices, and after each interaction, a measurement is done on $\mathcal R_p$. Then the above sequence $(X_p)_p$ represents the sequence of measurement outcomes, and the sequence $(\rho_p)_p$ represents the state of the physical system $\mathcal S$ after the first $p$ measurements (we refer the reader to section 6 of \cite{AGS})). Our results immediately give a law of large numbers, a central limit theorem and a large deviation principle for the statistics of the measurements $(X_p)_p$.

We will pay specific attention to the application of our results to the case where the internal state space of the particle (what is sometimes called the \emph{coin space}) is two-dimensional. This will allow us to illustrate the full structure of homogeneous OQRWs, and in particular the notions of irreducibility, period, as well as the Baumgartner-Narnhofer (\cite{BN}) decompositions discussed previously in~\cite{CP1}.

Of the above cited articles, some give a central limit theorem for the position~$(X_p)_p$ associated with an OQRW on~$\zz^d$. The most general result so far is given in \cite{AGS}, and its proof is based on a central limit theorem for martingales and the K\"ummerer-Maassen ergodic theorem (see \cite{KuMa}). Our proof is based on a completely different strategy, using a computation of the Laplace transform, and uses an irreducibility assumption which does not appear in existing central limit results. We will show, however, that the irreducibility assumption can be dropped in some situations, and that our central limit theorem contains the result of \cite{AGS}, but yields more general formulas.

In addition, we can prove a large deviation principle for the position process~$(X_p)_p$ associated with an homogeneous OQRW on a lattice. The technique we used, based on the application of the Perron-Frobenius theorem to a suitable deformed positive map, goes back (to the best of our knowledge) to \cite{HMO}. None of the articles cited above proves a large deviation principle. As we were completing this paper, however, we learnt of the recent article \cite{vHG}, which proves a similar result. We comment on this in section \ref{section_cltldp}.

The structure of the present paper is the following: in section \ref{section_OQRWs} we recall the main definitions of open quantum random walks 
%$\M$,
 specialized to the case where the underlying graph is a lattice in $\rr^d$, and define the auxiliary map of an open quantum random walk.
 %$\L$
In section~\ref{section_CPmaps} we recall standard results about irreducibility and period of completely positive maps. In section \ref{section_irreducibility} we characterize irreducibility and period of the open quantum random walk and its auxiliary map. In section \ref{section_cltldp} we state our main results: the central limit theorem and the large deviation principle. In section \ref{section_ZdC2} we specialize to the situation where the underlying graph is $\zz^d$ and the internal state space is $\cc^2$, and characterize each situation in terms of the transition operators.
%$L_s$ (see the next section for the precise definitions). 
In section \ref{section_examples} we study explicit examples, and compare our theoretical results to simulations.

\section{Homogeneous open quantum random walks}
\label{section_OQRWs}  

In this section we recall basic results and notations about open quantum random walks. We essentially follow the notation of \cite{CP1}, but specialize to the homogeneous case. For a more detailed exposition we refer the reader to \cite{APSS}.

We consider a Hilbert space $\mathfrak h$ and a locally finite lattice $V\subset \rr^d$, which we assume contains $0$, and is positively generated by a set $S\neq\{0\}$, in the sense that any $v$ in $V$ can be written as $s_1+\ldots+s_n$ with $s_1,\ldots, s_n \in S$. In particular, $V$ is an infinite subgroup of $\rr^d$. The canonical example is $V=\zz^d$, with $S=\{\pm v_1,\ldots, \pm v_d\}$ where $(v_1,\ldots,v_d)$ is the canonical basis of $\rr^d$.

We denote by $\mathcal H$ the Hilbert space $\mathcal H= \mathfrak h \otimes\cc^V$. We view $\mathcal H$ as describing the degrees of freedom of a particle constrained to move on $V$: the ``$V$-component" describes the spatial degrees of freedom (the position of the particle) while $\mathfrak h$ describes the internal degrees of freedom of the particle. According to quantum mechanical canon, we describe the state of the system as a positive, trace-class operator $\rho$ on $\H$ with trace one. Precisely, such an operator will be called a state.
\smallskip

We consider a map on the space $\mathcal I _1(\mathcal H)$ of trace-class operators, given by
\begin{equation}\label{eq_OQRW}
\mathfrak M \, : \, \rho \mapsto \sum_{j\in V}\, \sum_{s\in S}
\left(L_{s} \otimes \ketbra{j+s}{j}\right) \,\rho \,\left(L_{s}^* \otimes \ketbra{j}{j+s}\right)
\end{equation}
where the $L_{s}$, $s\in S$, are operators acting on $\h$ satisfying
\begin{equation}\label{eq_stochastic}
\sum_{s\in S} L_{s}^* \, L_{s}=\id.
\end{equation}
The $L_{s}$ are thought of as encoding both the probability of a transition by the vector $s$, and the effect of that transition on the internal degrees of freedom. Equation \eqref{eq_stochastic} therefore encodes the ``stochasticity" of the transitions.

\begin{remark}
The map $\M$ defined above is a special case of a quantum Markov chain, as introduced by Gudder in \cite{Gudder}. See Section 8 of \cite{CP1} for more comments.
\end{remark}

We associate with the OQRW $\M$ the auxiliary map $\L$ on the space $\mathcal I_1(\mathfrak h)$ of trace-class operators on $\h$ defined by  
\begin{equation}\label{eq_OQRWaux}
\mathfrak L \, : \,  \rho \mapsto\sum_{s\in S} L_{s}\,\rho \, L_{s}^*.
\end{equation}

Both \eqref{eq_OQRW} and \eqref{eq_OQRWaux} define trace-preserving (TP) maps, which are completely positive (CP), \textit{i.e.} for any $n$ in $\nn^*$, the extensions $\M\otimes \id$ and $\L\otimes\id$ to $\mathcal I_1(\mathcal H)\otimes \mathcal B(\cc^n)$ and $\mathcal I_1(\h)\otimes \mathcal B(\cc^n)$, respectively,  are positive. In particular, such a map transforms states (understood here as positive elements of $\mathcal I_1(\mathcal H)$ with trace one) into states. A completely-positive, trace-preserving map will be called a CP-TP map. We will call a map $\M$, as defined by \eqref{eq_OQRW}, an open quantum random walk, or OQRW; and we will call $\L$ the auxiliary map of~$\M$. 
To be more precise, we should call such an $\M$ an \emph{homogeneous} OQRW, but will drop the adjective homogeneous in the rest of this paper.

Let us recall that the topological dual $\mathcal I_1(\mathcal H)^*$ can be identified with $\mathcal B(\mathcal H)$ through the duality
\[(\rho,X)\mapsto \tr(\rho\, X).\]
\begin{remark}\label{remark_TPnormone}
When $\Phi=\M$ or $\Phi=\L$,
the adjoint $\Phi^*$ is a positive, unital (\textit{i.e.}~$\Phi^*(\id)=\id$) map on $\mathcal B (\mathcal H)$ (respectively $\mathcal B (\h)$), and by the Russo-Dye theorem (\cite{RD}) one has $\|\Phi^*\|=\|\Phi^*(\id)\|$ where the latter is the operator norm on $\mathcal B (\mathcal H)$ (respectively $\mathcal B (\mathcal H)$). This implies that trace-preserving positive maps have norm one, and in particular $\|\M\|=~1$ and $\|\L\|=1$.
\end{remark}

\begin{remark}\label{remark_classical} As noted in \cite{APSS}, classical Markov chains can be written as open quantum random walks. In the present case, if we have a subgroup $V$ of $\rr^d$ generated by a set $S$, and a Markov chain on $V$ with translation-invariant transition matrix $P=(p_{i,j})_{i,j\in V}$ induced by the law $(p_s)_{s\in S}$ on $S$, in the sense that $p_{i,j}=0$ if $j-i\not\in S$ and $p_{i,j}=p_{j-i}$ otherwise, then taking $\h=\cc$ and $L_{s}=\sqrt{p_{s}}$ 
induces the Markov chain with transition matrix $P$. This OQRW is called the minimal dilation of the Markov chain (see \cite{CP1} for a discussion of minimal and non-minimal dilations). Note that in this case the reduced map $\L$ is trivial:~$\L=1$.
\end{remark}

A crucial remark is that, for any initial state $\rho$ on $\mathcal H$, which is therefore of the form
\[ \rho=\sum_{i,j\in V} \rho(i,j)\otimes \ketbra ij ,\]
the evolved state $\mathfrak M(\rho)$ is of the form
\begin{equation}\label{eq_Mn}
\mathfrak M (\rho) = \sum_{i\in V} \mathfrak M (\rho,i) \otimes |i\rangle \langle i | , 
\ \mbox{ where }\ \M(\rho,i)=\sum_{s\in S} L_{s}\, \rho(i-s,i-s)\, L_{s}^*.
\end{equation}
Each $\mathfrak M(\rho,i)$ is a positive, trace-class operator on $\mathfrak h$ and
$ \sum_{i\in V} \tr\,\mathfrak M(\rho,i) =1$.
We notice that off-diagonal terms $\rho(i,j)$, for $i\neq j$, do not appear in $\M(\rho)$, and $\M(\rho)$ itself is diagonal. For this reason, from now on, we will only consider states of the form $\rho=\sum_{i\in V}\rho(i)\otimes \ketbra ii$. Equation \eqref{eq_Mn} remains valid, replacing $\rho(i,i)$ by~$\rho(i)$.
\smallskip

We now describe the (classical) processes of interest associated with $\M$. We begin with an informal discussion of these processes and their laws, and will only define the underlying probability space at the end of this section. We start from a state of the form $\rho=\sum_{i\in V}\rho(i) \otimes \ketbra ii$. We evolve $\rho$ for a time $p$, obtaining the state $\M^p(\rho)$ which, according to the previous discussion, is of the form
\[\M^p(\rho)=\sum_{i\in V} \M^p(\rho,i)\otimes\ketbra ii.\]
We then make a measurement of the position observable. According to standard rules of quantum measurement, we obtain the result $i\in V$ with probability $\tr\,\M^p(\rho,i)$. Therefore, the result of this measurement is a random variable $Q_p$, with law $\pp(Q_p=i)= \tr\,\mathfrak M^p(\rho,i)$ for $i\in V$. In addition, if the position $Q_p=i\in V$ is observed, then the state is transformed to $\frac{\mathfrak M^p(\rho,i)}{\tr \mathfrak M^p(\rho,i)}$. This process $(Q_p,\frac{\M^p(\rho,Q_p)}{\tr\,\M^p(\rho,Q_p)})$ we call the process ``without measurement", to emphasize the fact that virtually only one measurement is done, at time $p$. Notice that, in practice, two values of this process at times $p<p'$ cannot be considered simultaneously as the measure at time $p$ perturbs the system, and therefore subsequent measurements.

Now assume that we make a measurement at every time $p\in\nn$, applying the evolution by $\M$ between two measurements. Again assume that we start from a state $\rho$ of the form $\sum_{i\in V}\rho(i)\otimes |i\rangle \langle i|$. Suppose that at time $p$, the position was measured at $X_p=j$ and the state (after the measurement) is $\rho_p \otimes \ketbra {j}{j}$. Then, after the evolution, the state becomes
\[\M(\rho_p \otimes \ketbra {j}{j})= \sum_{s\in S} L_{s}\, \rho_p\, L_{s}^* \otimes \ketbra {j+s}{j+s},\]
so that a measurement at time $p+1$ gives a position $X_{p+1}=j+s$ with probability $\tr\, L_{s}\, \rho_p\, L_{s}^*$, and then the state becomes $\rho_{p+1}\otimes\ketbra {j+s}{j+s}$ with $\rho_{p+1}=\frac{L_{s}\, \rho_p\, L_{s}^*}{\tr\, L_{s}\, \rho_p\, L_{s}^*}$. The sequence of random variables $(X_p,\rho_p)$ is therefore a Markov process with transitions defined by
\begin{equation}\label{eq_Markov} 
\pp\Big((X_{p+1},\rho_{p+1})=(j+s,\frac{L_{s}\,\sigma\, L_{s}^*}{\tr(L_{s}\sigma L_{s}^*)})
\Big| (X_p,\,\rho_p)=(j,\sigma)\Big)=\tr (L_{s}\,\sigma\, L_{s}^*),\end{equation}
for any $j\in V$, $s\in S$ and $\sigma\in{\cal I}_1(\h)$ and initial law $\pp\big((X_0,\rho_0)=(i,\frac{\rho(i)}{\tr\rho(i)})\big)=\tr\rho(i).$
Note that the sequence $X_0=i_0$, \ldots, $X_p=i_p$ is observed with probability 
\begin{equation} \label{eq_probatraj}
\pp(X_0=i_0, \ldots, X_p=i_p)=\tr\,\big(L_{s_p}\ldots L_{s_1}L_{i_0}\,\rho(i_0)\, L_{i_0}^*L_{s_1}^*\ldots L_{s_p}^*\big)
\end{equation}
if $i_1-i_0=s_1$,\ldots $i_p-i_{p-1}=s_p$ belong to $S$, and zero otherwise. In addition, this sequence completely determines the state $\rho_p$:
\begin{equation}\label{eq_rhon} \rho_p = \frac{L_{s_p}\ldots L_{s_1}\,\rho(i_0)\, L_{s_1}^*\ldots L_{s_p}^*}{\tr\,L_{s_p}\ldots L_{s_0}\,\rho(i_0)\, L_{s_0}^*\ldots L_{s_p}^*}.
\end{equation}
As emphasized in \cite{APSS}, this implies that, for every $p$, the laws of $X_p$ and $Q_p$ are the same, \textit{i.e.} \[ \pp(X_p=i)=\pp(Q_p=i)\quad \forall i\in V.\]

We now construct a probability space to carry the processes just described. Fixing an open quantum random walk $\M$ on $V$ defined by operators $(L_s)_{s\in S}$ we define the set $\Omega=V^\nn$, equipped with the $\sigma$-field generated by cylinder sets. An element of $\Omega$ is denoted by $\omega=(\omega_k)_{k\in \nn}$ and we denote by $(X_p)_{p\in\nn}$ the coordinate maps.  For any state $\rho$ on $\H$ of the form $\rho=\sum_{i\in V}\rho(i)\otimes \ketbra ii$, we define a probability $\pp^{(p)}_\rho$ on $V^{p+1}$ by formula \eqref{eq_probatraj}. One easily shows, using the stochasticity property \eqref{eq_stochastic}, that the family $(\pp_\rho^{(p)})_p$ is consistent, and can therefore be extended uniquely to a probability $\pp_\rho$ on $\Omega$.  We denote by $\rho_p$ the random variable
\[ \rho_p =  \frac{L_{X_p- X_{p-1}}\ldots L_{X_1-X_{0}} \,\rho(X_0)\,L^*_{X_1- X_{0}}  \ldots L^*_{X_p- X_{p-1}}}{\tr(L_{X_p- X_{p-1}}\ldots L_{X_1- X_{0}} \,\rho(X_0)\,L^*_{X_1-X_{0}}  \ldots L^*_{X_p- X_{p-1}})}.\] 
We will also denote $Q_p=X_p$, but will only use the notation $Q_p$ when we consider ``non-measurement" experiments, and in particular will never consider an event implying simultaneously outcomes $Q_p$ and $Q_{p'}$ for $p\neq p'$. These processes reproduce the behaviour of the measurement outcomes and of the associated resulting states. In particular, equation \eqref{eq_Markov} above holds in a mathematical sense with $\pp_\rho$ replacing $\pp$. From now on, we will usually drop the $\rho$ in $\pp_\rho$.

\section{Irreducibility and period: general results}
\label{section_CPmaps}

In this section we focus on the general notions of irreducibility and period for a completely positive (CP) map $\Phi$ on $\mathcal I_1(\mathcal K)$, where $\mathcal K$ is a separable Hilbert space which, in practice, will be either $\h$ or $\H$. We assume $\Phi$ is given in the form
\begin{equation}\label{eq_defKraus}
\Phi(\rho)=\sum_{\kappa\in K} A_\kappa \rho A_\kappa^*
\end{equation}
where $K$ is a countable set, and the series $\sum_{\kappa\in K} A_\kappa^* A_\kappa$ is strongly convergent. This is the case for operators such as $\M$ or $\L$  and we actually know from the Kraus theorem that this is the case for any completely positive $\Phi$, see \cite{Kraus} or \cite{NieChu}, where this is called the operator-sum representation.
We recall that such a map is automatically bounded as a linear map on $\mathcal I_1(\mathcal K)$ (see \textit{e.g.} Lemma~2.2 in \cite{Sch}), so that it is also weak-continuous. In most practical cases, we will additionally assume that $\|\Phi\|=1$; this will be the case, in particular, if~$\Phi$ is trace-preserving.
\smallskip

We give various equivalent definitions of the notion of irreducibility for $\Phi$ which was originally defined by Davies in \cite{Dav}. Note that this original definition holds for $\Phi$ positive, but for simplicity, we discuss it only for maps $\Phi$ which are completely positive (CP), and therefore have a Kraus decomposition \eqref{eq_defKraus}. The equivalence between the different definitions, as well as the relevant references, are discussed in \cite{CP1}.  We recall some standard notations: an operator $X$ on $\mathcal K$ is called positive, denoted $X\geq0$, if, for $\phi\in \mathcal K$,  one has $\langle \phi, X\, \phi\rangle \geq 0$. It is called strictly positive, denoted $X>0$, if, for $\phi\in \mathcal K\setminus\{0\}$, one has $\langle\phi, X\, \phi\rangle>0$.
\begin{defi}\label{defi_irreducibility}
The CP map $\Phi$ is called irreducible if one of the following equivalent conditions hold:
\begin{itemize}
\item the only orthogonal projections $P$ reducing $\Phi$, \textrm{i.e.} such that~$\Phi\big(P\mathcal I_1(\mathcal K)P\big) \subset P\mathcal I_1(\mathcal K)P$, are $P=0$ and $\id$,
\item for any $\rho\geq 0$, $\rho\neq 0$ in $\mathcal I_1(\mathcal K)$, there exists $t$ such that~$\e^{t \Phi} (\rho)>~0$,
\item for any non-zero $\phi\in\mathcal K$, the set
$ \cc[A] \, \phi $ is dense in $\mathcal K$, where $\cc[A]$ is the set of polynomials in $A_\kappa, \, \kappa\in K$,
\item the only subspaces of $\mathcal K$ that are invariant by all operators $A_\kappa$ are $\{0\}$ and~$\mathcal K$.
\end{itemize}
\end{defi}

\begin{remark}
Note also that the notion of irreducibility is strongly related to the notion of subharmonic projection, again see \cite{CP1}.
\end{remark}
We will also need on occasion the notion of regularity, which is evidently stronger than irreducibility:
\begin{defi}\label{defi_regularity}
The CP map $\Phi$ is called $N$-regular if one of the equivalent conditions hold:
\begin{itemize}
\item for any $\rho\geq 0$, $\rho\neq 0$ in $\mathcal I_1(\mathcal K)$, one has $\Phi^N (\rho)>~0$,
\item for any non-zero $\phi\in\mathcal K$, the set
$\{ A_{\kappa_1}\ldots A_{\kappa_N}\, \phi\,|\, \kappa_1,\ldots,\kappa_N\in K \}$ is total in $\mathcal K$.
\end{itemize}
The map $\Phi$ is called regular if it is $N$-regular for some $N$ in $\nn^*$. 
\end{defi}

\begin{remark}\label{remark_regularity}
The following properties are quite immediate:
\begin{itemize} 
\item If $\Phi$ is irreducible, then $\vee_{\kappa\in K} {\rm Ran} A_\kappa={\mathcal K}$ (while the converse is not true).
\item If $\vee_{\kappa\in K} {\rm Ran} A_\kappa={\mathcal K}$ and $\sigma$ is a faithful state, then $\Phi(\sigma)$ is faithful.
Indeed, we can write $\sigma=\sum_j\sigma_j |u_j\rangle\langle u_j|$, with $\sigma_j>0$ and $(u_j)_j$ an orthonormal basis for $\mathcal K$.
Then $\Phi(\sigma)=\sum_{j,\kappa}\sigma_j |A_\kappa u_j\rangle\langle A_\kappa u_j|$, and the conclusion easily follows.
\item If $\vee_{\kappa\in K} {\rm Ran} A_\kappa={\mathcal K}$ and $\Phi$ in $N$-regular, $N\ge 1$, then $\Phi$ is $(N+n)$\,--\,regular for any $n\ge 0$.
This is an immediate consequence of the previous point.
\end{itemize}
\end{remark}

The following proposition, which is a Perron-Frobenius theorem for positive maps on $\mathcal I_1(\mathcal K)$, essentially comes from \cite{EHK} (for the finite dimensional case) and~\cite{Sch} (for the infinite dimensional case). To state it in sufficient generality, we need to recall the definition of the spectral radius of a map $\Phi$:
\[r(\Phi)=\sup \{|\lambda|, \, \lambda\in \mathrm{Sp}\,\Phi\}\]
where $\mathrm{Sp}\,\Phi$ is the spectrum of $\Phi$.
\begin{prop}\label{prop_Schrader}
Assume a CP map $\Phi$ on $\mathcal I_1(\mathcal K)$ has an eigenvalue $\lambda$ of modulus $r(\Phi)$, with eigenvector $\rho$, and
either $\mathrm{dim}\, \mathcal K <\infty$ or $r(\Phi)=\|\Phi\|$.
Then:
\begin{itemize}
\item $|\lambda|$ is also an eigenvalue, with eigenvector $|\rho|$,
\item if $\Phi$ is irreducible, then $\mathrm{dim}\,\mathrm{Ker}\,(\Phi-\lambda\,\id)=1$.
\end{itemize}
In particular, if $\Phi$ is irreducible and has an eigenvalue of modulus $r(\Phi)$,
then~$r(\Phi)$ is an eigenvalue with geometric multiplicity one, with an eigenvector that is a strictly positive operator.
\end{prop}
\begin{remark}\label{rem_finitedim}
When $\Phi$ is a completely positive, trace-preserving map, one has $\|\Phi\|=1$, so that the conclusion applies if $\lambda$ is of modulus $1$. In \cite{CP1}, this was enough, since we applied this result to the operator $\M$. In section \ref{section_cltldp} we will also need to apply it to a deformation of the operator $\L$, which will no longer be trace-preserving.
\end{remark}

\begin{remark}
The previous proposition gives in particular uniqueness and faithfulness of the invariant state, when it exists, for an irreducible map $\Phi$.
As one can expect, the converse result holds: if $\Phi$ admits a unique invariant state and that state is faithful, then $\Phi$ is irreducible (see \cite{CP1}, Section 7).
\end{remark}

\medskip
We now turn to the notion of period for positive maps. We will denote by~$\submod$,~$\addmod$ the substraction and addition \emph{modulo} $d$.
\begin{defi}\label{def_period}
Let $\Phi$ be a CP, trace-preserving, irreducible map and consider a resolution of the identity
$(P_0,\ldots,P_{d-1})$, \emph{i.e.} a family of orthogonal
projections such that $\sum_{j=0}^{d-1} P_j=\id$. One says that $(P_0,\ldots,P_{d-1})$
is $\Phi$-cyclic if $P_j A_\kappa= A_\kappa P_{j\submod 1}$ for $j=0,\ldots,d-1$ and any $k$. 
The supremum of all $d$ for which there exists a $\Phi$-cyclic resolution of identity $(P_0,\ldots,P_{d-1})$ is called the \textrm{period} of $\Phi$.
If $\Phi$ has period~$1$ then we call it aperiodic.
\end{defi}
\begin{remark}
If $\mathrm{dim}\,\mathcal K$ is finite then the period is always finite.
\end{remark}

The following proposition is the analog of a standard result for classical Markov chains:
\begin{prop}\label{prop-periodic-blocks}
Assume $\Phi$ is completely positive, irreducible, with finite period $d$, and denote by $P_0,\ldots,P_{d-1}$ a cyclic decomposition of $\Phi$. Then :
\begin{enumerate}
\item we have the relation $
\Phi(P_i\, \rho \, P_j)=P_{i\addmod 1}\,  \Phi(\rho) \, P_{j\addmod 1},
$
\item 
for any $j=0,\ldots, d-1$, the restriction $\Phi^d_j$ of $\Phi^d$ to $P_j{\cal I}_1({\mathcal K})P_j$ is irreducible aperiodic,
\item if $\Phi$ has an invariant state $\rhoinv$, then $\Phi^d_j$ has a unique invariant state 
$\rhoinv_j\overset{\mathrm{def}}=d\times P_j\rhoinv P_j$.
\end{enumerate}\end{prop}

\pre

\begin{enumerate}
\item The first relation is obvious, and shows that $P_j{\cal I}_1({\mathcal K})P_j$ is stable by $\Phi^d$.
\item Consider a state $P_j\rho P_j$ in $P_j{\cal I}_1({\mathcal K})P_j$. By irreducibility of $\Phi$, $\e^{t\Phi}(P_j\rho P_j)$ is faithful, so~$P_j\e^{t\Phi}(P_j \rho P_j)P_j$ is faithful in $\mathrm{Ran}\,P_j$. But by the relation in point $1$,
\[P_j\e^{t\Phi}(P_j\rho P_j)P_j =\sum_{n=0}^\infty  \frac{t^{dn}}{(dn)!}\, \Phi^{dn}(P_j \rho P_j) = \sum_{n=0}^\infty  \frac{t^{dn}}{(dn)!}\, (\Phi_j^{d})^n(P_j \rho P_j).\]
This shows that $\Phi_j^d$ is irreducible. Now, if $\Phi_j^d$ has a cyclic decomposition of identity $(P_{j,0},\ldots, P_{j,\delta-1})$ then by the commutation relations this induces a cyclic decomposition of identity for $\Phi$ with $d\times \delta$ elements. Therefore,~$\delta=1$.
\item The invariance of $\rhoinv_j$ is trivial by point 1, and the irreducibility of $\Phi_j^d$ implies the unicity of the invariant state. By remark 4.8 in \cite{CP1}, $\tr(P_j \rhoinv P_j)$ does not depend on $j$, so it is $1/d$.
\fin
\end{enumerate}
The following results were originally proved by Fagnola and Pellicer in \cite{FP} (with partial results going back to \cite{EHK} and \cite{Gro}). We recall that the point spectrum of an operator is its set of eigenvalues, and that we denote by $\mathrm{Sp}_{pp}\Phi^*$ the point spectrum of $\Phi^*$.
\begin{prop}\label{prop_aperiodicCPTP}
If $\Phi$ is an irreducible, completely positive, trace-preserving map on~$\mathcal I_1(\mathcal K)$ and has finite period $d$ then:
\begin{itemize}
\item the set $\mathrm{Sp}_{pp}\Phi^* $, is a subgroup of the circle group $\mathbb T$,
\item the primitive root of unity $\e^{\i 2\pi/d}$ belongs to $\mathrm{Sp}_{pp}\Phi^*$ if and only if $\Phi$ is~$d$-periodic.
\end{itemize}
\end{prop}
An immediate consequence is the following:
\begin{prop}\label{prop_cvgCPTPaperiodic}
If a completely positive, trace-preserving map $\Phi$ on $\mathcal I_1(\mathcal K)$ is irreducible and aperiodic with invariant state $\rhoinv$, and $\mathcal K$ is finite-dimensional then
\begin{itemize}
\item $\mathrm{Sp_{pp}}\,\Phi\cap\mathbb T = \{1\}$,
\item for any $\rho\in \mathcal I_1(\mathcal K)$ one has $\Phi^p(\rho)\rightarrow \rhoinv$ as $p\to\infty$.
\end{itemize}
\end{prop}

\section{Irreducibility and period of $\M$ and $\L$}\label{section_irreducibility}

Now we turn to the case where the operator $\Phi$ is an open quantum random walk~$\M$ generated by $L_s$, $s\in S$, or the auxiliary map $\L$ as defined by \eqref{eq_OQRWaux}.
We will study irreducibility and periodicity properties of both operators $\M$ and $\L$ and mutual relations. We will first explain why we focus on a study of $\L$, when~$\M$ should intuitively be the object of interest.

For any $v$ in $V$ we denote
\[\mathcal P_\ell(v)= \{\pi=(s_1,\ldots, s_\ell)\in S^{\ell}\, |\, \sum_{p=1}^\ell s_p = v\}\]
and, in addition, we consider  
\[\mathcal P(v)=\cup_{\ell\geq 1}\mathcal P_\ell(v),
\qquad
\mathcal   P_\ell = \cup_{v\in V} \mathcal P_\ell(v) \qquad \mathcal P=\cup_{\ell\in \nn}\mathcal P_\ell = \cup_{v\in V} \mathcal P(v).\]
In analogy with \cite{CP1}, we use the notation
\[
L_\pi = L_{s_\ell}\cdots L_{s_1},
\qquad
\mbox{ for } \pi=(s_1,\ldots, s_\ell)\in \mathcal P_\ell.
\]
We remark that the notations for the paths and the set of paths are slightly different from our previous paper \cite{CP1} since we can use homogeneity, which allows us to drop the dependence on the particular starting point.

The irreducibility of $\L$ and $\M$ are easily characterized in terms of paths. This is true in general for OQRWs (see \cite{CP1}, Proposition 3.9 in particular), but the following characterization for $\M$ is specific to homogeneous OQRWs.
\begin{lemma}\label{lemma_irreducLM}
Let $\M$ be an open quantum random walks defined by transition operators $L_s$, $s\in S$, and $\L$ its auxiliary map.
\begin{enumerate}
\item The operator $\L$ is irreducible if and only if, for any $x\neq 0$ in $\h$, the set~$\{L_\pi x,\, |\, \pi \in \mathcal P\}$ is total in~$\h$.
\item
The operator $\M$ is irreducible if and only if, for any $x\neq 0$ in $\h$ and $v$ in~$V$, the set $\{L_\pi x,\, |\, \pi \in \mathcal P(v)\}$ is total in $\h$.
\end{enumerate}
\end{lemma}
\pre

This lemma is proven by a direct application of Definition \ref{defi_irreducibility} (third condition). For $\L$, it is immediate. For $\M$, one sees easily that irreducibility amounts to the fact that, for any $x\otimes \ket {w}$, the set $\{L_\pi x \otimes \ket{v+w}, \pi\in\mathcal P(v)\}$ is dense in~$\h$ for any $v\in V$ (see the details in Proposition 3.9 of \cite{CP1}) and this is equivalent to the statement above.
\fin
\smallskip

This lemma obviously implies the following result:
\begin{coro}\label{coro_ML}
If $\M$ is irreducible, then $\L$ is irreducible.
\end{coro}
One can, however, prove a more explicit criterion for irreducibility of $\M$ than Lemma \ref{lemma_irreducLM}. For consistency we rephrase here one of the equivalent definitions in Definition \ref{defi_irreducibility}:
\begin{prop}\label{prop_caracM}
The operator $\L$ is irreducible if and only if the operators
\[\{L_{s},\ s\in S\}\]
have no invariant closed subspace in common, apart from $\{0\}$ and $\h$.

The operator $\M$ is irreducible if and only if the operators
\[\{L_{\pi_0},\ \pi_0\in \mathcal P(0)\}\]
have no invariant closed subspace in common, apart from $\{0\}$ and $\h$.
\end{prop}

\pre

The characterization for $\L$ immediately follows from the last condition in Definition \ref{defi_irreducibility}. If $\M$ is not irreducible then for some $v\in V$, the closed space 
\[ \h_v=\overline{\mathrm{Vect}\{L_\pi x,\, \pi \in \mathcal P(v)\}} \]
is different from $\h$. Since the concatenation of any $\pi_0\in \mathcal P(0)$ with $\pi\in\mathcal P(v)$ gives an element of $\mathcal P(v)$, the space $\h_v$ must be $L_{\pi_0}$-invariant. Conversely, if all operators $L_{\pi_0}$ have an invariant subspace $\mathfrak h'$ in common, then for any $x\in \mathfrak h'$, the set $\{L_\pi x,\, |\, \pi \in \mathcal P(0)\}$ is contained in $\mathfrak h'$ and $\M$ is not irreducible.  \fin

Proposition \ref{prop_caracM} allows us to construct examples of OQRWs such that $\L$ is irreducible, but not $\M$.
\begin{example}\label{ex_LirrMnotirr}
Let $d=1$ and 
\[L_+=\begin{pmatrix}0 & a_+ \\ b_+ & 0 \end{pmatrix}\qquad L_-=\begin{pmatrix}0 & a_- \\ b_- & 0 \end{pmatrix}\]
with $a_+, a_-, b_+, b_-$ positive, with $a_+^2+b_+^2=a_-^2+b_-^2=1$ and $a_+ b_- \neq a_- b_+$. Then, by Proposition \ref{prop_caracM}, 
$\L$ is irreducible, but $\M$ is not, since the vectors of the canonical basis are eigenvectors for any $L_{\pi}$, $\pi\in\mathcal P(0)$ (see also Proposition \ref{reducibility-C2}). 
\end{example}

The following proposition is proved in \cite{CP1}. We reprove it here.
\begin{prop}\label{prop_inexistenceetatinv}
Assume $\M$ is irreducible. Then it does not have an invariant state.
\end{prop}

\pre

By Corollary \ref{coro_ML}, $\L$ is irreducible, so it has a unique invariant state $\rhoinv$ on~$\h$, which is faithful. Assume $\M$ has an invariant state; by irreducibility it is unique. Since $\M$ is translation-invariant, any translation of that state would be also invariant, so by unicity the invariant state is translation-invariant. It must then be of the form $\sum_{i\in V}\rhoinv\otimes \ketbra ii$, but this has infinite trace, a contradiction.
\fin

All ergodic convergence results for $\M$ given in \cite{CP1} assume the existence of an invariant state.
This is similar to the situation for classical Markov chains; however, some interesting asymptotic properties  of $\M$ can be studied in the absence of an invariant state, and this includes large deviations or central limit theorems. As we will see, such properties  can be derived from the study of $\L$. This is why, in the study of homogeneous OQRWs, the focus shifts from $\M$ to~$\L$.
%We recall here that, as we wrote in Remark \ref{remark_classical}, this map $\L$ is trivial for minimal dilations of classical Markov chains.
\smallskip

To avoid discussing trivial cases, in the rest of this paper we will usually make the following assumption, which by Remark \ref{remark_regularity} automatically holds as soon as $\L$ (or $\M$) is irreducible:
\begin{center}
\textbf{Assumption H1:} one has the equality $ \overline{\bigvee_{s\in S} \mathrm{Ran}\,L_{s}}=\h$.
\end{center}
\smallskip

This assumption is a natural one, since after just one step, even in the reducible case, the system is effectively restricted to the space $\bigvee_{s\in S} \mathrm{Ran}\,L_{s} $. More precisely, for any positive operator $\rho$ on $\h$, one has, for any $s$, 
\[\mathrm{supp}\, L_s\, \rho\, L_s^* \subset \mathrm{supp}\, \L(\rho)\subset \overline{\bigvee_{s\in S} \mathrm{Ran}\,L_{s}}.\]

Note that we have not given results equivalent to Lemma \ref{lemma_irreducLM} for the notion of regularity. We do this here:
\begin{lemma}\label{lemma_regularLM}
The operator $\L$ is $N$-regular if and only if for any $x\neq 0$ in $\h$, the set $\{L_\pi x,\, |\, \pi \in \mathcal P_N\}$ is total in~$\h$.
The operator $\M$ can never be regular.
\end{lemma}

\pre

This is obtained by direct application of Definition \ref{defi_regularity}, that shows the criterion for $\L$. It also shows that $\M$ is $N$-regular if and only if for any $x\neq 0$ in $\h$, any $v$ in $V$, the set $\{L_\pi x,\, |\, \pi \in \mathcal P_N(v)\}$ is total in~$\h$. However, if the distance from the origin to $v$ is larger than $N$, then $\mathcal P_N(v)$ is empty. \fin

\smallskip
One could be tempted to consider a weaker version of regularity for $\L$ where the index $N$ can depend on $\rho$. The following result shows that, if $\h$ is finite-dimensional, this is not weaker than regularity:
\begin{lemma}\label{lemma_regularity}
Assume $\h$ is finite-dimensional. If for every $\rho\geq 0$ in $\mathcal I_1(\h)\setminus\{0\}$, there exists $N>0$ such that $\L^N(\rho)$ is faithful, then there exists $N_0$ such that $\L$ is $N_0$-regular.
\end{lemma}

\pre

First observe that $\L$ is necessarily irreducible and so assumption H1 must hold.
%, since for any $\rho$ in~$\mathcal I_1(\h)$, $\mathrm{supp}\, \L(\rho)\subset \overline{\bigvee_{s\in S} \mathrm{Ran}\,L_{s}}$.
Besides, the current assumption implies that, for any $x$ in $\h$, there exists $N_x>0$ such that $\L^{N_x}(\ketbra xx)$ is faithful. Since faithfulness of $\L^{N_x}(\ketbra xx)$ is equivalent to the existence of a family $\pi_1,\ldots, \pi_{\mathrm{dim}\,\h}$ of paths of length $N_x$, such that the determinant of $(L_{\pi_1} x, \ldots, L_{\pi_{\mathrm{dim}\,\h}} x)$ is nonzero, there exist open subsets $B_x$ of the unit ball, such that $x\in B_{x_0}$ implies that $\L^{N_{x_0}}(\ketbra xx)$ is faithful. By compactness of the unit ball, there exists a finite covering by $B_{x_1}\cup\ldots\cup B_{x_p}$. Remark \ref{remark_regularity} then implies that if we let~$N_0=\sup_{i=1,\ldots,p}N_{x_i}$ one has $\L^{N_0}(\ketbra xx)$ faithful for any nonzero $x$. This implies that $\L$ is $N_0$-regular.
\fin
\medskip

We now turn to the notion of period for $\L$ and $\M$. By Definition \ref{def_period}, a resolution  of identity $(p_0, \ldots, p_{d-1})$ of $\h$ will be $\L$-cyclic if and only if
\begin{equation*}
p_j L_s = L_s p_{j\submod 1}\quad \mbox{for }j=0,\ldots,d-1 \mbox{ and any }s\in S.
\end{equation*}
Consequently, by Proposition \ref{prop-periodic-blocks}, we have
\begin{equation} \label{eq_commutationLu}
\L(p_j \, \rho \, p_j) =  p_{j\addmod 1}\, \L(\rho) \, p_{j\addmod 1}.
\end{equation}

\begin{remark}\label{remark_dsmallerdimh}
Since the $p_j$ sum up to $\id_\h$, the period of $\L$ cannot be greater than $\mathrm{dim}\,\h$, a feature which will be extremely useful when $\mathrm{dim}\,\h$ is small.
\end{remark}
On the other hand, as we observed in \cite{CP1}, a resolution  of identity $(P_0, \ldots, P_{d-1})$ of $\H$ will be $\M$-cyclic if and only if it is of the form
\begin{equation}\label{eq_cyclicity}
P_k=\sum_{i\in V}P_{k,i}\otimes \ketbra ii \quad\mbox{with}\quad
P_{k,i} L_s = L_{s} P_{k\submod 1,i+s}.
\end{equation}

\begin{remark}
The cyclic resolutions for $\M$ are translation invariant, in the sense that, if $P_k=\sum_{i\in V}P_{k,i}\otimes \ketbra ii $, $k=0,...d-1$, is a cyclic resolution for $\M$, then also $P'_k=\sum_{i\in V}P_{k,i+v}\otimes \ketbra ii $, $k=0,...d-1$, is a cyclic resolution for any $v$.
\end{remark}

We will, however, make little use for a cyclic resolution of identity for $\M$ in this paper. On the other hand, the periodicity of $\L$ can be an easy source of information on $\M$:

\begin{prop}\label{prop-M-period}
We have the following properties:
\begin{enumerate}
\item The period of $\M$, when finite, is even.
\item If $\L$ is irreducible and has even period $d$, then $\M$ is reducible.
\end{enumerate}
\end{prop}

\pre

\begin{enumerate}
\item Assume that $(P_0,\ldots,P_{d-1})$ is a $\M$-cyclic resolution of identity associated with $\M$. As we observed above, the $P_k$ are of the form
\[P_k=\sum_{i\in V}P_{k,i}\otimes \ketbra ii \quad\mbox{with}\quad
P_{k,i} L_s = L_{s} P_{k\submod 1,i+s}.\]
Then if we call $i$ in $V$ odd or even depending on the parity of its distance to the origin, define 
\[P_{k,\textrm{odd}}=\sum_{i\, \textrm{odd}}P_{k,i}\otimes \ketbra ii \quad \mbox{and} \quad P_{k,\textrm{even}}=\sum_{i\, \textrm{even}}P_{k,i}.\]
Then $(P_{0,\textrm{odd}}, P_{1,\textrm{even}}, P_{2,\textrm{odd}},\ldots)$ is a $\M$-cyclic resolution of identity.
\item Denote by $(p_0,\ldots,p_{d-1})$ a cyclic resolution of identity associated with $\L$. Define 
\[p_{\mathrm{odd}}=\sum_{k\,\mathrm{odd}} p_k\qquad p_{\mathrm{even}}=\sum_{k\,\mathrm{even}} p_k.\]
It is obvious from relations \eqref{eq_cyclicity} that $\mathrm{Ran}\, p_{\mathrm{odd}}$ and $\mathrm{Ran}\, p_{\mathrm{even}}$ are nontrivial invariant spaces for any $L_{\pi_0}$, $\pi_0\in\mathcal P(0)$. We conclude by Proposition \ref{prop_caracM}. \fin
\end{enumerate}

\smallskip

Last, we give an analogue of a classical property of Markov chains with finite state space:
\begin{lemma}\label{lemma_irrapereg}
If $\h$ is finite-dimensional, then the map $\L$ is irreducible and aperiodic if and only if it is regular.
\end{lemma}

\pre

If $\L$ is irreducible and aperiodic, then by Proposition \ref{prop_cvgCPTPaperiodic} for any state $\rho$ on $\h$, one has $\L^p(\rho)\underset{n\to\infty} {\longrightarrow}\rhoinv$ so that $\L^p(\rho)$ is faithful for large enough $p$. By Lemma \ref{lemma_regularity}, this implies the regularity of $\L$. Conversely, if $\L$ is regular, then it is irreducible, and for any projection $p$, the operator $\L^N(p)$ is faithful, so that $p$ cannot be a member of a cyclic resolution of identity unless $p=\id$.
\fin

\section{Central Limit Theorem and Large Deviations }\label{section_cltldp}

The Perron-Frobenius theorem for CP maps allows us to obtain a large deviations principle and a central limit theorem for the position process $(X_p)_{p\in\nn}$ (or, equivalently, for the process $(Q_p)_{p\in\nn}$) associated with an open quantum random walk $\M$ and an initial state $\rho$ (see section \ref{section_OQRWs}). In most of our statements, we assume for simplicity that $\L$ is irreducible. We discuss extensions of our results at the end of this section.

Before going into the details of the proof, we should mention that, as we were completing the present article, we learnt about the recent paper \cite{vHG}, which proves a  large deviation result for empirical measures of outputs of quantum Markov chains, which can be viewed as the ``steps" $(X_p-X_{p-1})_p$ taken by an open quantum random walk. This result is similar to the statement in our Remark \ref{remark_Sanov}, and implies a level-1 large deviation result for the position $(X_p)_p$ when the OQRW is irreducible and aperiodic. In addition, the statement in \cite{vHG} extends to a large deviations principle for empirical measures of $m$-tuples of $(X_p-X_{p-1})_p$. Our (independent) result, however, treats the case where the OQRW is irreducible but not aperiodic, and can be extended beyond the irreducible case.

For the proofs of this section, it will be convenient to introduce some new notations. For $u$ in $\rr^d$ we define $L_{s}^{(u)}= \e^{\braket us /2}L_s$, and denote $\L_u$ the map induced by the $L_{s}^{(u)}$, $s\in S$: for $\rho$ in $\mathcal I_1(\h)$,
\[\L_u(\rho)=\sum_s L_{s}^{(u)}\rho L_{s}^{(u)*}.
\]

The operators $\L_u$ will be useful in order to treat the moment generating functions of the random variables $(X_p)$:
\begin{lemma}\label{def_Lu}
For any $u$ in $\rr^d$ one has 
\begin{equation}\label{mgf}
\ee(\exp\,\braket u {X_p-X_0}) = \sum_{i_0\in V}  \tr\big(\L_{u}^p(\rho(i_0))\big).
\end{equation}
\end{lemma}

\pre

For any $k$ in $\nn^*$ let $S_k=X_{k+1}-X_k$ and consider $u\in\rr^d$. Then we have
\begin{eqnarray*}
&& \hspace{-2em}\ee(\exp\,\braket u{X_p-X_0})\\
&=& \sum_{i_0\in V}\sum_{s_1,\ldots,s_p\in S^p} \hspace{-1em}\pp(X_0=i_0,S_1=s_1,\ldots, S_p=s_p)\, \exp\,\braket{u}{s_1+\ldots+s_p} \\
&=& \sum_{i_0\in V}\sum_{s_1,\ldots,s_p\in S^p} \hspace{-1em}\tr(L_{s_p}\ldots L_{s_1}\, \rho(i_0)\, L_{s_1}^*\ldots L_{s_p}^*)\, \exp\,\braket{u}{s_1+\ldots+s_p} \end{eqnarray*}
and this gives formula \eqref{mgf}.
\fin
\begin{remark}\label{remark_XpQp}
One also has
\[ \ee(\exp\,\braket u{X_p})= \ee(\exp\,\braket u{Q_p})=\sum_{i_0\in V} \exp\braket{u}{i_0}\,\tr\big(\L_{u}^p(\rho(i_0))\big).\]
This will allow us to give results analogous to Theorem \ref{theo_ldp} and \ref{theo_clt} for the process~$(Q_p)_p$. Note that considering $X_p$ or $X_p-X_0$ is essentially equivalent, but as we remarked in section \ref{section_OQRWs}, $Q_p$ and $Q_0$ cannot be considered simultaneously.
\end{remark}
\smallskip

The following lemma describes the properties of the largest eigenvalue of $\L_u$:
\begin{lemma}\label{lemma_analyticitylambdau}
Assume that $\h$ is finite-dimensional and $\L$ is irreducible.  For any $u$ in $\rr$, the spectral radius $\lambda_u\overset{\mathrm{def}}=r(\L_u)$ of $\L_u$ is an algebraically simple eigenvalue of $\L_u$, and has an eigenvector $\rho_u$ which is a strictly positive operator, and we can normalize it to be a state. In addition, the map $u\mapsto \lambda_u$ can be extended to be analytic in a neighbourhood of~$\rr^d$.
\end{lemma}
\smallskip

\pre

By Lemma \ref{lemma_irreducLM}, if $\L$ is irreducible, then so is any $\L_u$ for $u\in \rr^d$. Proposition~\ref{prop_Schrader}, applied here specifically to an Hilbert space of finite dimension, gives the first sentence except for the algebraic simplicity of the eigenvector $\lambda_u$, as it implies only the geometric simplicity.
If we can prove that, for all $u$ in $\rr^d$, the eigenvalue $\lambda_u$ is actually algebraically simple  then the theory of perturbation of matrix eigenvalues (see Chapter II in \cite{Kato}) will give us the second sentence. 
Now, in order to prove the missing point, consider the adjoint $\L_u^*$ of $\L_u$ on $\mathcal B(\h)$, which in this finite-dimensional setting, can be identified, with $\mathcal I_1(\h)$. It is easy to see from Definition \ref{defi_irreducibility} that $\L_u^*$ is irreducible. Its largest eigenvalue is $\lambda_u$, with eigenvector $M_u$, which, by Proposition \ref{prop_Schrader}, is invertible. We can consider the map 
$$
\widetilde \L _u : \rho\mapsto \frac1{\lambda_u}\,M_u^{1/2} \L_u( M_u^{-1/2}\, \rho \, M_u^{-1/2}) M_u^{1/2}.
$$
This $\widetilde \L _u$ is clearly completely positive, and is trace-preserving since $\widetilde \L _u^*(\id)=\id$.
%map on $\mathcal I_1(\h)$, for the preservation of the trace, we can directly verify that, for $\rho\in \mathcal I_1(\h)$,
%$$
%\tr(\widetilde \L _u(\rho))
%= \frac1{\lambda_u}  \tr(\L_u^*(M_u)( M_u^{-1/2}\, \rho \, M_u^{-1/2}) )
%= \tr(\rho).
%$$
Proposition \ref{prop_Schrader} shows that $\widetilde \L _u$ has $1$ as a geometrically simple eigenvalue, with a strictly positive eigenvector $\widetilde \rho _u$. Then $1$ must also be algebraically simple, otherwise there exists $\eta_u$ such that $\widetilde \L _u(\eta_u)= \eta_u + \widetilde\rho_u$, but taking the trace of this equality yields $\tr(\widetilde\rho _u)=0$, a contradiction. This implies that $\L_u$ has $\lambda_u$  as a algebraically simple eigenvalue.
\fin
\smallskip

We can now state our large deviation result:
\begin{theo}\label{theo_ldp}
Assume that $\h$ is finite-dimensional and that $\L$ is irreducible. Then the process $(\frac1p(X_p-X_0))_{p\in\nn^*}$ associated with $\M$ satisfies a large deviation principle with a good rate function~$I$. Explicitly, there exists a lower semicontinuous mapping $I:\mathbb R^d\to\mathbb [0, + \infty]$ with compact level sets $\{x\, |\, I(x)\leq \alpha\}$, such that, for any open $G$ and closed~$F$ with $G\subset F \subset \mathbb R^d$, one has
\begin{eqnarray*}
&&-\inf_{x\in G} I(x) \leq \liminf_{p\to\infty}\frac1p \log P(\frac {X_p-X_0}{p}\in G)
\\&& \qquad\qquad\qquad
\leq\limsup_{p\to\infty}\frac1p \log P(\frac {X_p-X_0}{p}\in F)\leq -\inf_{x\in F} I(x).
\end{eqnarray*}
\end{theo}
\begin{remark}\label{remark_LDPXp}
If we add the assumption that $X_0$ has an everywhere defined moment generating function, {\it e.g.} that the initial state $\rho$ satisfies 
$\ee(\exp\,\braket u{X_0})= \sum_{i_0\in V} \e^{\braket u{i_0}}\tr\rho(i_0)< \infty$ for all $u$ in $\rr^d$, then this theorem also holds for $(X_p)_p$ or equivalently $(Q_p)_p$ in place of $(X_p-X_0)_p$.
\end{remark}

\begin{remark}\label{remark_GvH}
Using the techniques detailed in \cite{vHG}, it is possible, for any $m$ in $\nn$, to extend the above theorem and obtain a full large deviation principle for the sequence of $(m+1)$-tuples $\frac1p(X_p-X_0, X_{p+1}-X_1,\ldots, X_{p+m}-X_m)_p$, or (under the same condition as in Remark \ref{remark_LDPXp}) for $(X_p,\ldots,X_{p+m})_p$.
\end{remark}

\pre

We start with equation \eqref{mgf}. Since $\mathfrak h$ is finite-dimensional, if $\rho(i_0)$ is faithful, then, with $r_{u,i_0}=\inf\mathrm{Sp}(\rho(i_0))>0$ and $s_{u,i_0}=\frac{\tr\rho(i_0)}{\inf\mathrm{Sp}(\rho_u)}>0$ (where $\mathrm{Sp}(\sigma)$ denotes the spectrum of an operator $\sigma$),
\begin{equation}\label{eq_preuveldp1}
r_{u,i_0}\, \rho_u \leq \rho(i_0) \leq s_{u,i_0} \,\rho_u.
\end{equation}
Note that $r_{u,i_0}\leq \tr \rho(i_0)$ so that both $r_{u,i_0}$ and  $s_{u,i_0}$ are summable along $i_0$. Consequently, we shall have
\begin{equation}\label{eq_preuveldp2}
r_{u,i_0} \, \lambda_u^p\,\rho_u   \leq \L_u^p\big(\rho(i_0) \big)\leq s_{u,i_0} \,\lambda_u^p\, \rho_u .
\end{equation}
Using these bounds in relation \eqref{mgf}, we immediately obtain, for all $u\in\rr^d$,
\[
 \lambda_u^p\,\sum_{i_0\in V} r_{u,i_0}\, \rho_u
\le \ee(\exp\,\braket u {X_p-X_0}) \le   \lambda_u^p\,\sum_{i_0\in V} s_{u,i_0}\, \rho_u
\]
where the sums are finite and strictly positive; so that
\begin{equation}\label{eq_GE}\lim_{p\to\infty}\frac1p \, \log\ee(\exp\,\braket u {X_p})  = \log\lambda_u .\end{equation}

Now, if $\rho(i_0)$ is not faithful, but $\L$ is aperiodic, due to Proposition \ref{prop_cvgCPTPaperiodic}, then $\L^N(\rho(i_0))$ is faithful for large enough $N$, and \eqref{eq_preuveldp1} holds with $\L_u^N(\rho(i_0))$ in place of $\rho(i_0)$ and \eqref{eq_preuveldp2} holds with $(p-N)$ instead of $p$ in the exponents of $\lambda_u$. We still recover \eqref{eq_GE}.

Finally, if $\rho(i_0)$ is not faithful and $\L$ has period $d>1$, then, considering a cyclic decomposition of identity $(p_0,\ldots,p_{d-1})$, we can consider the single blocks of the form $p_j\rho(i_0)p_j$.
By Proposition \ref{prop-periodic-blocks}, $\L^d$ is irreducible aperiodic when restricted to each $p_j{\cal I}_1({\h})p_j$ 
and $\L_u^d(p_j \rho_u p_j) = \lambda_u^d p_{j} \rho_u p_{j}$.
Then, by the regularity of the restrictions of $\L^d$, using Remark \ref{remark_regularity} and the obvious extension of \eqref{eq_commutationLu} to $\L_u$, there exist $N\in {\mathbb N}$ and $r_{u,i_0},s_{u,i_0}>0$ such that, for any block $p_j\rho(i_0)p_j \neq 0$,
\begin{equation*}%\label{eq_preuveldp1-2}
r_{u,i_0}\, p_j\,\rho_u \,p_j \leq p_j\, \L_u^{dN}\rho(i_0)\,p_j \leq s_{u,i_0} \,p_j\, \rho_u\, p_j
\end{equation*}
and if $p=dN+r$, $r\in\{0,\ldots,d-1\} $,
\begin{equation*}%\label{eq_preuveldp2-2}
r_{u,i_0}\, \lambda_u^{p-dN}\,p_{j+r}\,\rho_u \,p_{j+r} \leq \L_u^p\big(p_j \,\rho(i_0)\,p_j \big)
\leq s_{u,i_0} \,\lambda_u^{p-dN}\, p_{j+r}\,\rho_u \,p_{j+r}.
\end{equation*}
Summing over $j$, we recover equation \eqref{eq_GE} again.

In any case, we obtain \eqref{eq_GE} for all $u\in\rr^d$. Lemma \ref{lemma_analyticitylambdau} shows that $u\mapsto \log\lambda_u$ is analytic on $\rr$. Applying the G\"artner-Ellis theorem (see \cite{DZ}) we obtain the bounds mentioned in Theorem \ref{theo_ldp}, with rate function
\[I(x) = \sup_{u\in\rr^d} \, \big(\langle u,x \rangle - \log\lambda_u\big).\ {\Box}\] 

\begin{remark}\label{remark_Sanov}
If $\varphi$ is any function $S\to \rr$ and $S_p= \sum_{k=1}^p \varphi(X_k-X_{k-1})$ then the process~$(\frac{S_p}p)_{p\in\nn}$ also satisfies a large deviation principle, with rate function
\[I_\varphi(x) = \sup_{t\in\rr} \, \big(t\, x  - \log\lambda_{t \varphi}\big)\] 
where  $\lambda_{t\varphi}$ is the largest eigenvalue of 
\[\L_{t \varphi} : \rho \mapsto \sum_{s\in S}\e^{ t \varphi(s)} L_s \, \rho L_s^*.  \]
This is shown by an immediate extension of the proofs of Lemma \ref{lemma_analyticitylambdau} and Theorem \ref{theo_ldp}, and yields a level-2 large deviation result for the process $(X_p-X_{p-1})$ (\emph{e.g.} using Kifer's theorem \cite{Kifer}).
\end{remark}

\begin{remark}
As noted in Remark \ref{remark_classical}, when $\M$ is the minimal dilation of a classical Markov chain with transition probabilities $(p_s)_{s\in S}$, the map $\L$ is trivial: it is just multiplication by $1$ on $\rr$. The maps $\L_u$, however, are not trivial: they are multiplication by
\[\lambda_u = \sum_{s\in S}\exp\braket us \, p_s.\]
We therefore recover the same rate function as in the classical case, see \textit{e.g.} section 3.1.1 of \cite{DZ}. 
\end{remark}

\begin{remark}
The technique of applying the Perron-Frobenius theorem to a $u$-dependent deformation of the completely positive map defining the dynamics, goes back (to the best of our knowledge) to \cite{HMO}, and is a non-commutative adaptation of a standard proof for Markov chains.
\end{remark}

We denote by $c$ the map $c:\rr^d\ni u\mapsto \log\lambda_u$. As is well-known (see \textit{e.g.} section II.6 in \cite{Ellis}), the differentiability of $c$ at zero is related to a law of large numbers for the process $(X_p)_{p\in \nn}$. Similarly, the second order differential will be relevant for the central limit theorem. 

\begin{coro}\label{coro_derivatives}
Assume that $\h$ has finite dimension and that $\L$ is irreducible. The function $c$ on $\rr^d$ is infinitely differentiable at zero. Denote by 
\[\L_u' : \rho \mapsto \sum_{s\in S} \braket us \, L_s \rho L_s^* \quad \mbox{ and }\quad\L_u'' : \rho \mapsto \sum_{s\in S} {\braket us}^2  L_s \rho L_s^*.\]
Then, denoting $\lambda_u'\overset{\mathrm{def}}{=}{\frac\d{\d t}}_{|t=0} \lambda_{tu}$ and $\lambda_u''\overset{\mathrm{def}}{=}\frac{\d^2}{\d t^2}_{|t=0} \lambda_{tu}$, we have
\begin{equation}\label{eq_lambdaprime}
\lambda_u' =\tr\big(\L_u'(\rhoinv)\big)
\end{equation}
\begin{equation}\label{eq_lambdasec}
\lambda_u'' = \tr\big(\L''_u(\rhoinv)\big)+2\tr\big(\L'_u(\eta_u)\big)
\qquad
\end{equation}
where $\eta_u$ is the unique solution with trace zero of the equation
\begin{equation}\label{eq_etau}
\big(\id-\L\big)(\eta_u) = \L'_u(\rhoinv)-\tr\big(\L'_u(\rhoinv)\big)\,\rhoinv.
\end{equation}
This implies immediately that
\begin{equation}\label{eq_diffc}
\d c (0)\,(u)=\lambda_u' \qquad  \d^2 c(0)\, (u,u)=\lambda_u''-\lambda_u'{}^2.
\end{equation}
\end{coro}

\pre

Lemma \ref{lemma_analyticitylambdau} shows that $c_u$ is infinitely differentiable at any $u\in \rr^d$.
In addition (again see Chapter II in \cite{Kato}), the largest eigenvalue  $\lambda_u$ of $\L_u$ is an analytic perturbation of $\lambda_0=1$, and has an eigenvector $\rho_{u}$ which we can choose to be a state, and this $\rho_u$ is an analytic perturbation of $\rho_0$. Then one has
\[\lambda_{tu} = 1 + t \lambda'_u + \frac {t^2}2 \lambda_u'' + o(t^2)\]
\[ \rho_{tu} = \rhoinv + t \,\eta_u + \frac {t^2}2 \,\sigma_u+o(t^2)\]
\[ \L_{tu} = \L + t \,\L'_u + \frac{t^2}2\, \L''_u + o(t^2)\]
and since every $\rho_{tu}$ is a state then $\tr \,\eta_u=\tr\, \sigma_u=0$. Then the relation $\L_{tu}(\rho_{tu})=\lambda_{tu}\, \rho_{tu}$ 
yields 
\[\L'_u(\rhoinv) + \L(\eta_u) = \eta_u + \lambda'_u\, \rhoinv\]
\[\frac12\L(\sigma_u)+ \L_u'(\eta_u)+ \frac12 \L_u''(\rhoinv)= \frac12 \sigma_u + \lambda_u'\, \eta_u + \frac12 \lambda_u '' \, \rhoinv. \]
Taking the trace of the first relation immediately yields relation \eqref{eq_lambdaprime}. In addition, it yields relation \eqref{eq_etau}. Since $\id - \L$ has kernel of dimension one, and range in the set of operators with zero trace, it induces a bijection on that state, so that \eqref{eq_etau} has a unique solution with trace zero. Then taking the trace of the second relation above, and using the fact that $\L$ is trace-preserving gives relation \eqref{eq_lambdasec}.
\fin

\begin{coro}\label{coro_lln}
Assume that $\h$ has finite dimension and $\L$ is irreducible, and let~$m=\sum_{s} \tr(L_s\rhoinv L_s^*) \, s$. Then the process $(\frac1p(X_p-X_0))_{p\in\nn}$  associated with $\M$ converges exponentially to $m$, \emph{i.e.} for any $\eps>0$ there exists $N>0$ such that, for large enough $p$, 
\[\pp(\big\|\frac{X_p-X_0}p-m\big\|>\eps)\leq \exp -pN.\]
This implies the almost-sure convergence of $(\frac{X_p}p)_{p\in\nn}$ to $m$.
\end{coro}
\begin{remark}
The almost-sure convergence holds replacing~$X_p$ by $Q_p$.
\end{remark}

\pre

This is a standard result, see \emph{e.g.} Theorem II.6.3 and Theorem II.6.4 in \cite{Ellis}.
\fin

\begin{theo}\label{theo_clt}
Assume that $\h$ is finite-dimensional and $\L$ is irreducible. Denote by $m$ the quantity defined in Corollary \ref{coro_lln}, and by $C$ the covariance matrix associated with the quadratic form $u\mapsto \lambda_u'' - \lambda_u'{}^2$.
 Then the position process~$(X_p)_{p\in\nn}$ associated with $\M$ satisfies
\[ \frac{X_p-p\,m}{\sqrt p} \underset{p\to\infty}\longrightarrow \mathcal N(0,C)\]
where convergence is in law.
\end{theo}

\begin{remark}
Again this result holds replacing $X_p$ by $Q_p$.
\end{remark}
\begin{remark}
The formulas for the mean and variance are the same as in \cite{AGS} when $V=\zz^d$ and $S=\{\pm v_i,\, i=1,\ldots,d\}$ ($v_1,\ldots,v_d$ is the canonical basis of $\rr^d$). This can be observed from the fact that, if $Y_u$ is the unique (up to a constant multiple of the $\id$) solution of equation 
\[(\id-\L^*)(Y_u)= \sum_{s\in S}\braket us L_s^*L_s - \braket um \,\id,\]
(note that our $Y_u$ is the $L_l$ of \cite{AGS}) then
\[\tr\big(\L_u'(\eta_u)\big) = \tr\big(\L_u'(\rhoinv) Y_u\big) - \tr\big(\L_u'(\rhoinv)\big)\,\tr\big(\rhoinv Y_u\big)\]
and denoting $Y_i=Y_{v_i}$ we have
\begin{eqnarray*}
\braket u {Cu} = \sum_{i,j=1}^d u_i u_j &&\hspace{-2em}\Big(\ind_{i=j}\big(\tr (L_{+i}\rhoinv L_{+i}^*) +\tr (L_{-i}\rhoinv L_{-i}^*) \big)\\
&& +2\tr (L_{+i}\rhoinv L_{+i}^*\,Y_j) - 2\tr (L_{-i}\rhoinv L_{-i}^*\, Y_j)\\
&& - 2m_i \tr(\rhoinv \, Y_j)- m_im_j\Big)
\end{eqnarray*}
which leads to the formula for $C$ given in \cite{AGS}:
\begin{eqnarray*}
C_{i,j}= &&\hspace{-2em}\ind_{i=j}\big(\tr (L_{+i}\rhoinv L_{+i}^*) +\tr (L_{-i}\rhoinv L_{-i}^*) \big)\\
&&\hspace{-2em} +\big(\tr (L_{+i}\rhoinv L_{+i}^* \,Y_j) +\tr (L_{+j}\rhoinv L_{+j}^*\,Y_i)\big)\\
&&\hspace{-2em} \big( \tr (L_{-i}\rhoinv L_{-i}^*\, Y_j)+\tr (L_{-j}\rhoinv L_{-j}^*\, Y_i)\big)\\
&&\hspace{-2em} - \big(m_i \tr(\rhoinv \, Y_j)+m_j \tr(\rhoinv \, Y_i) \big)
- m_im_j.
\end{eqnarray*}
\end{remark}

%\begin{remark}
%It is possible to derive the central limit result from our study of the function $u\mapsto \log \ee(\exp \braket u {X_p})$ and a multivariate version of Bryc's theorem (see Appendix A.4 in \cite{QSM}). This is, however not immediate because, even though it is clear that $\ee(\exp \braket u {X_p})$ can be extended analytically to a complex neighbourhood of the origin, it is not so obvious that this is possible for the function $\log \ee(\exp \braket u {X_p})$. We choose, however, to consider a proof of Theorem \ref{theo_clt} that does not use Bryc's theorem because it leads to extensions beyond the case of irreducible $\L$. These extensions are discussed at the end of this section. 
%\end{remark}

\noindent \textbf{Proof of Theorem \ref{theo_clt}:}

Let us first consider the case where $\L$ is irreducible and aperiodic.  Equation \eqref{mgf} implies
\[  \ee(\exp\braket u {X_p-X_0})=\sum_{i_0\in V} \tr\big(\L_u^p(\rho(i_0))\big). \]
Now, considering the Jordan form of $\L$ shows that, if \[\delta\overset{\mathrm{def}}=\sup\{|\lambda|,\lambda\in\mathrm{Sp}\,\L\setminus\{1\}\},\]
then $\delta<1$ and for $u$ in a real neighbourhood of $0$ and $p$ in $\nn$,
\begin{equation}\label{eq_Jordan} \L_u^p = \lambda_u^p  \big(\varphi_u(\cdot)\, \rho_u + O((\delta+\varepsilon)^p)\big) \end{equation}
for some $\varepsilon$ such that $\delta+\varepsilon <1$, where $\varphi_u$ is a linear form on $\mathcal I_1(\h)$, analytic in $u$ and such that $\varphi_0=\tr$ and the $O((\delta+\varepsilon)^p)$ is in terms of the operator norm on $\mathcal I_1(\h)$. This implies
\begin{equation}\label{eq_mgf2}
\frac1p\log \sum_{i_0\in V}\tr (\L_u^p(\rho(i_0))) = \log \lambda_u + \frac1p \log \sum_{i_0\in V}\varphi_u(\rho(i_0)) + O((\delta+\eps)^p)
\end{equation}
for $u$ in the above real neighbourhood of the origin. This and Lemma \ref{lemma_analyticitylambdau} implies that the identity
\begin{equation}\label{eq_mgf3}
\lim_{p\to\infty}\frac 1p\log\ee(\exp\braket u {X_p-X_0}) = \log \lambda_u
\end{equation}
holds for $u$ in a neighbourhood of the origin. In addition, by  equation \eqref{eq_mgf2} and Corollary \ref{coro_derivatives},
\[ \lim_{p\to\infty}\frac 1{ p} \big(\nabla \log\ee(\exp\braket u {X_p-X_0})-pm)=0 \qquad \lim_{p\to\infty}\frac 1p \nabla^2 \log\ee(\exp\braket u {X_p-X_0})=C.\]
By an application of the multivariate version of Bryc's theorem (see Appendix A.4 in \cite{QSM}), we deduce that 
\[ \frac{X_p-X_0-p\,m}{\sqrt p} \underset{p\to\infty}\longrightarrow \mathcal N(0,C)\]
and this proves our statement in the case where $\L$ is irreducible aperiodic.

We now consider the case where $\L$ is irreducible with period $d$. Let $p_0,\ldots, p_{d-1}$ be a cyclic partition of identity; then, writing $p=qd+r$ we have for any $i_0 \in V$
\begin{eqnarray*}
\tr\big(\L_u^p(\rho(i_0))\big)&=& \sum_{j=0}^{d-1}\tr\big(p_j\,\L_u^{qd+r}(\rho(i_0)) \, p_j\big)\\
&=& \sum_{j=0}^{d-1}\tr\big(\L_u^{qd}(p_j\,\L_u^r(\rho(i_0)) \, p_j)\big)
\end{eqnarray*}
by a straightforward extension of \eqref{eq_cyclicity} to $\L_u$. By Proposition \ref{prop-periodic-blocks}, for any $j$, $r$ and the previous discussion, one has
\[\lim_{q\to\infty}\frac1{qd}  \log \tr\big(\L_u^{qd}(p_j\,\L_u^r(\rho(i_0)) \, p_j)\big) = \log \lambda_u\]
and one can extend all terms in this identity so that it holds in a complex neighbourhood of the origin. This finishes the proof of our statement. \fin

\begin{remark}
The reader might wonder why we need to go through the trouble of considering relations \eqref{eq_Jordan} and \eqref{eq_mgf2} to derive the extension of  \eqref{eq_mgf3} to complex $u$. This is because there is no determination of the complex logarithm that allows to consider $\log\ee(\exp\braket u {X_p-X_0})$ for complex $u$ and arbitrarily large $p$. This forces us to start by transforming $\frac1p\log\ee(\exp\braket u {X_p-X_0})$.
\end{remark}

\paragraph{Generalizations of  Theorems \ref{theo_ldp} and \ref{theo_clt}}

We finish with a discussion of possible generalizations of Theorems \ref{theo_ldp} and \ref{theo_clt} beyond the case of irreducible $\L$. To this aim, we introduce the following subspaces of $\h$:
\begin{equation}\label{eq_defDR}
\mathcal D = \{\phi\in \h\, |\,\langle\phi,\L^p(\rho)\,\phi\rangle\underset{p\to\infty}{\longrightarrow}0 \mbox{ for any state }\rho\} \quad \mbox{and} \quad \R=\D^\perp.
\end{equation}
Alternatively, $\R$ can be defined as the supremum of the supports of $\L$-invariant states, and $\D$ as $\R^\perp$. Note in particular that $\mathrm{dim}\,\R\geq 1$ and $\R$ is invariant by all operators $L_s$, $s\in S$. These subspaces are the Baumgartner-Narnhofer decomposition of $\h$ associated with $\L$ (see \cite{BN} or \cite{CP1}).
Note that, in \cite{CP1}, we only considered the spaces $\D_\M$ and $\R_\M$  associated with $\M$ instead of $\L$. Here the decomposition for $\M$ plays no role and $\R_\M$ is equal to $\{0\}$.

The following result will replace the Perron-Frobenius theorem when $\L$ is not irreducible. The proof can be easily adapted from Proposition \ref{prop_Schrader} and Lemma~\ref{lemma_analyticitylambdau}.
\begin{prop}\label{prop_PFnonirr}
The following properties are equivalent:
\begin{enumerate}
\item the auxiliary map $\L$ has a unique invariant state $\rhoinv$,
\item the restriction $\L_{|\mathcal I_1(\R)}$ of $\L$ to $\R$ is irreducible,
\item the value $1$ is an eigenvalue of $\L$ with algebraic multiplicity one.
\end{enumerate}
If, in addition, $\L_{|\mathcal I_1(\R)}$ is aperiodic, then $1$ is the only eigenvalue of modulus one, and for any state $\rho$, one has $\L^p(\rho)\underset{p\to\infty}\longrightarrow \rhoinv$.
\end{prop}

This leads to an extension of Theorem \ref{theo_clt} to the cases
\begin{itemize}
\item where $\L_{|\mathcal I_1(\mathcal R)}$ is irreducible (even if $\R \neq\h$); by Proposition \ref{prop_PFnonirr}, this is equivalent to $\L$ having a unique invariant state;
\item when $\mathcal R=\h$.% has a unique decomposition as $\bigoplus \R_k $ where each $\L_{|\mathcal I_1(\R_k)}$ is irreducible.
\end{itemize}
With these two extensions, our central limit theorem has the same generality as the one given in \cite{AGS}: the first case is Theorem 5.2 of that reference, the second case is treated in Section $7$ in \cite{AGS}.
These extensions are proven observing that:
\begin{itemize}
\item by Proposition \ref{prop_PFnonirr}, the proof of Theorem \ref{theo_clt} can immediately be extended to the situation where $\L_{|\mathcal I_1(\R)}$ is irreducible aperiodic, and from there to the situation where $\L_{|\mathcal I_1(\R)}$ is irreducible periodic;
\item when $\mathcal R=\h$, it admits a decomposition $\R=\oplus_k \R_k$ (see \cite{CP1}), each $\mathcal I_1(\R_k)$ is stable by $\L$, the restrictions $\L_{|\mathcal I_1(\R_k)}$ are irreducible, and the non-diagonal blocks do not appear in a probability like \eqref{eq_probatraj}.
\end{itemize}
We have seen in \cite{CP1} that one can always decompose $\h$ into $\h=\D\oplus \bigoplus_{k\in K} \R_k$
with each $\R_k$ as discussed above. However, in the general case, we do not have a clear statement of Theorems \ref{theo_ldp} and \ref{theo_clt} because 
if $\D$ is non-trivial and $\mathrm{card}\,K\geq 2$, it is difficult to control how the mass of $\rho_0$ will flow from $\D$ into the different components $\R_k$. 

Last, remark that the proof of Theorem \ref{theo_ldp} relies on the fact that $\L_u$ is irreducible. This holds if $\L$ is irreducible; the converse, however, is not true, and $\R_\L$ may be different from $\R_{\L_u}$. The proof of Theorem \ref{theo_ldp} can be extended to  derive a lower large deviation bound in the case when $\R=\h$ using the idea described above, but when $\L$ is not irreducible, the quantity $\lambda_u$ may not be analytic, in which case we \emph{a priori} obtain only the upper large deviation bound, see Example \ref{ex_ldpbreakdown}. 

\section{Open quantum random walks with lattice $\zz^d$ and internal space $\cc^2$}
\label{section_ZdC2}
The goal of this section is to illustrate our various concepts, and give explicit formulas in the case where $V=\zz^d$ and $\h=\cc^2$. We start with a study of the operators $\L$ and $\M$, and a characterization of their (ir)reducibility and of the associated decompositions of the state space in this specific situation.

We begin in Proposition \ref{prop_caracL} with a classification of the possible situations depending on the dimension of $\R$ (as defined in \eqref{eq_defDR}) and its possible decompositions. Then, in Lemma \ref{lemma_diag} we characterize those situations in terms of the form of the operators $L_s$. Later on, we also consider the period.
To avoid discussing trivial cases, we will make a second assumption:
\begin{center}\textbf{Assumption H2:} the operators $L_s$ are not all proportional to the identity.\end{center}
This is equivalent to saying that we assume $\L\neq \id$.
\medskip

We start by discussing the possible forms of $\R$ and $\D$:
\begin{prop}\label{prop_caracL}
Consider the operators $L_s$, $s\in S$, defining the open quantum random walk $\M$, and suppose that assumptions \emph{H1} and \emph{H2} hold. Then we are in one of the following three situations.
\begin{enumerate}
\item If the $L_{s}$ have no eigenvector in common, then $\L$ is irreducible, there exists a unique $\L$-invariant state which is faithful, and one has
\[\R = \h \qquad \D=\{0\}.\]
\item If the $L_{s}$ have only one (up to multiplication) eigenvector $e_1$ in common, then $\L$ is not irreducible, the state $\ketbra {e_1}{e_1}$ (if $\|e_1\|=1$) is the unique~$\L$-invariant state, and for any nonzero vector $e_2\perp e_1$, one has
\[ \R = \cc\, e_1\qquad \D = \cc\, e_2.\] 
\item If the $L_{s}$ have two linearly independent eigenvectors $e_1$ and $e_2$ in common, any invariant state is of the form $\rhoinv=t\,\ketbra {e_1}{e_1} + (1-t)\ketbra {e_2}{e_2}$ for $t\in[0,1]$, and one has
\[\R = \h = \cc\, e_1 \oplus \cc\, e_2 \qquad \D=\{0\}.\]
\end{enumerate}
\end{prop}

\pre

We recall that, by the fourth equivalent statement in Definition \ref{defi_irreducibility}, the map~$\L$ is irreducible if and only if the $L_{s}$ do not have a common, nontrivial, invariant subspace. If~$\h=\cc^2$ then this is equivalent to saying that the $L_{s}$ do not have a common eigenvector.

Now assume that $\L$ is not irreducible, so that the $L_{s}$ have a common norm one eigenvector $e_1$, with $L_{s}\,e_1=\alpha_{s}\, e_1$ for all $s$. Then $\ketbra {e_1}{e_1}$ is an invariant state. Complete $(e_1)$ into an orthonormal basis $(e_1,e_2)$. Then, if $\rho$ is an invariant state, $\rho=\sum_{i,j=1,2}\rho_{i,j}\,\ketbra{e_i}{e_j}$, and
\[ \L(\rho)= \sum_{i,j=1,2}\, \sum_{s\in S}\rho_{i,j}\,\ketbra{L_s e_i}{L_s e_j}.\]
Then 
\[\rho_{2,2}=\braket{e_2}{\rho\, e_2}=\sum_{s\in S}\rho_{2,2}\,|\braket{e_2}{L_s e_2}|^2\]
so that either $\rho_{2,2}=0$ or $\sum_{s\in S}|\braket{e_2}{L_s e_2}|^2=1$; but, since $\sum_{s\in S}\|L_s e_2\|^2=1$, this is possible only if $e_2$ is an eigenvector of all $L_s$, $s\in S$.

Therefore, in situation 2, $\ketbra{e_1}{e_1}$ is the only invariant state. In situation~3, observe that if there existed an invariant state with $\rho_{1,2}=\overline{\rho_{2,1}}\neq 0$, then any state would be invariant and $\L$ would be the identity operator, a case we excluded.
\fin

\begin{remark}
In situations $2$ and $3$ we recover the fact, proven in \cite{CP1} (and originally in \cite{BN}) that, if $\ketbra{e_1}{e_1}$ is an invariant state and $e_2\neq 0$ is in $e_1^\perp\cap\R$ then $\ketbra{e_2}{e_2}$ is an invariant state. The above proposition gives an explicit Baumgartner-Narnhofer decomposition of $\h$ (see \cite{BN} or sections 6 and 7 of \cite{CP1}).
In the case where $\h=\cc^2$, it turns out that $\R$ can always be written in a unique way as $\R=\bigoplus \R_k$ with $\L_{|\mathcal I_1(\R_k)}$ irreducible (except for the trivial case when $\L$ is the identity map). This is not true in general and is a peculiarity related to the low dimension of $\h$.
\end{remark}
\medskip

Next we study the explicit form of the operators $L_s$ in each of the situations described by Proposition \ref{prop_caracL}. We will use the standard notation that, for  two families of scalars $(\alpha_s)_{s\in S}$ and $(\beta_s)_{s\in S}$, $\|\alpha\|^2$ is $\sum_{s\in S}|\alpha_s|^2$ and $\braket \alpha \beta$ is~$\sum_{s\in S}\overline{\alpha_s}\beta_s$.  
\begin{lemma}\label{lemma_diag}
With the assumptions and notations of Proposition \ref{prop_caracL}:
\begin{itemize}
\item We are in situation $2$ if and only if there exists an orthonormal basis of~$\h=\cc^2$ in which
\[L_{s}=\begin{pmatrix}\alpha_{s} & \gamma_{s} \\ 0 & \beta_{s}\end{pmatrix}\]
for every $s$ with 
\[\|\alpha\|^2 =\|\beta\|^2+\|\gamma\|^2=1, \qquad \braket\alpha\gamma=0,\]
\[\sup_{s\in S} |\beta_{s}| >0,\qquad \sup_{s\in S} |\gamma_{s}| >0,\]
\[\mbox{ there exist } s\neq s'\mbox{ in } S \mbox{ such that }(\alpha_s-\beta_s)\,\gamma_{s'}\neq (\alpha_{s'}-\beta_{s'})\,\gamma_{s}.\]

\item We are in situation $3$ if and only if there exists an orthonormal basis of~$\h=\cc^2$ in which
\[L_{s}=\begin{pmatrix}\alpha_{s} &0 \\ 0 & \beta_{s}\end{pmatrix}\]
for every $s$, with 
\[\|\alpha\|^2 =\|\beta\|^2=1,\]
\[\mbox{there exists }s\mbox{ in } S \mbox{ such that } \alpha_s \neq \beta_s.\]\end{itemize}
\end{lemma}

\pre

This is immediate by examination.
\fin

\begin{remark}
In situation 2, let $\rho$ be any state. One has
\[ \braket {e_2} {\L^p(\rho)\, e_2} = \tr\big(\rho\, \L^{*\,p}(\ketbra {e_2}{e_2})\big) = \|\beta\|^{2p}\, \braket {e_2} {\rho\, e_2} \underset{p\to\infty}{\longrightarrow} 0  \]
by the observation that $  \|\beta\|^2 <1$. We recover the fact that $\D=\cc\, e_2$.
\end{remark}
\medskip

We now turn to the study of periodicity for the operator $\L$. We start with a simple remark:
\begin{remark}\label{rem_irrC2}
Whenever the operators $L_s$ have a common eigenvector $e$, then the restriction of $\L$ to $\mathcal I_1(\cc e)$ is aperiodic. In particular, if $\L$ is not irreducible but has a unique invariant state, then by necessity $\R$ is one-dimensional so that $\L_{|\R}$ must be aperiodic.
\end{remark}

In more generality, because $\mathrm{dim}\,\h=2$, by Remark \ref{remark_dsmallerdimh}, any irreducible~$\L$ has period either one or two. The following lemma characterizes those $L_s$ defining an operator $\L$ with period 2:
\begin{lemma}\label{lemma_Laperiodic}
The map $\L$ is irreducible periodic if and only if there exists a basis of $\h$ for which every operator $L_{s}$ is of the form $\begin{pmatrix} 0 & \gamma_{s}\\ \nu_{s}& 0\end{pmatrix}$. In that case, for any~$s\neq s'$, one has $\gamma_s \,\nu_{s'}\neq \gamma_{s'}\,\nu_s$ and $\|\gamma\|^2=\|\nu\|^2=1$, and the unique invariant state of $\L$ is $\frac12\,\id$.
\end{lemma}

\pre

If the period of $\L$ is two, then the cyclic partition of identity must be of the form $\ketbra {e_1}{e_1}, \ketbra {e_2}{e_2}$ and the cyclicity imposes the relations
\[L_{s} e_1 \in \cc\, e_2, \quad L_{s} e_2 \in \cc\, e_1 \quad \mbox{for any }s\in S.\]
This gives the form of the $L_s$. The condition $\sum_s|\gamma_s|^2=\sum_s|\nu_s|^2=1$ simply follows by the trace preservation property. Now observe that the eigenvalues of $L_s$ are solutions of $\lambda_s^2=\gamma_s \nu_s$. Fix one solution $\lambda_s$, the other being $-\lambda_s$. Then $\begin{pmatrix}x\\y\end{pmatrix}$ is an eigenvector if and only if $\gamma_s y=\pm \lambda_s \, x$.
Therefore, two operators $L_s$ and $L_{s'}$ have an eigenvector in common if and only if $\nu_s\,\lambda_{s'}=\pm\nu_{s'}\,\lambda_{s}$. This is easily seen to be equivalent to $\gamma_s \,\nu_{s'}= \gamma_{s'}\,\nu_s$. Last, one easily sees that the equation
\[\sum_{s\in S}L_s \begin{pmatrix}a&b\\c&d\end{pmatrix} L_s^* = \begin{pmatrix}a&b\\c&d\end{pmatrix}\]
is equivalent to $a=d$, $b=\braket \nu \gamma\, c$ and  $c=\braket  \gamma\nu \, b$. Moreover, $|\braket  \gamma\nu|=1$ would imply that the vectors $(\gamma_s)_{s\in S}$ and $(\nu_s)_{s\in S}$ are proportional, which is forbidden by irreducibility. Therefore $a=d$ and $b=c=0$.
\fin

\medskip
The following theorem is a central limit theorem for all open quantum random walks satisfying H1 and H2. It gives more explicit expressions for the parameters of the limiting Gaussian, except when $\L$ is irreducible aperiodic, in which case the parameters of the Gaussian are given in Theorem~\ref{theo_clt}.
\begin{theo}\label{theo_hisc2}
Assume an open quantum random walk with $V=\zz^d$ and $\h=\cc^2$ satisfies assumptions \emph{H1}, \emph{H2}. Then there exist $m\in\cc^d$ and $C$ a $d\times d$ positive semi-definite matrix such that we have the convergence in law
\[ \frac{X_p-p\,m}{\sqrt p} \underset{p\to\infty}\longrightarrow \mathcal N(0,C).\]
Following the notation of Lemmas \ref{lemma_diag} and \ref{lemma_Laperiodic} we have:
\begin{itemize}
\item In situation 1, if $\L$ is periodic, consider two random variables $A$ and~$B$ with $\pp(A=s)=|\nu_s|^2$ and $\pp(B=s)=|\gamma_s|^2$. Then we have 
\[m=\frac12(\ee(A)+\ee(B))\qquad C=\frac12(\mathrm{var}(A)+\mathrm{var}(B)).\]
\item In situation 2, consider a classical random variable $A$ with $\pp(A=s)=~|\alpha_s|^2.$ Then we have
\[m=\ee(A) \qquad C=\mathrm{var}(A).\]
\item In situation 3, consider two classical random variables $A$ and~$B$ with $\pp(A=s)=|\alpha_s|^2$ and $\pp(B=s)=|\beta_s|^2$, and denote 
$p=\sum_{i\in V}\braket{e_1}{\rho(i)\,e_1}$, where $\rho$ is the initial state. Then we have
\[m=p\,\ee(A)+(1-p)\, \ee(B) \qquad C=p\,\mathrm{var}(A)+(1-p)\, \mathrm{var}(B).\] 
\end{itemize}
\end{theo}

\pre 

If $\L$ is irreducible periodic, for any $\sigma=\begin{pmatrix}\sigma_{11} & \sigma_{12} \\ \sigma_{21} & \sigma_{22} \end{pmatrix}$, we have
$$
L_s \sigma L_s^* = 
\begin{pmatrix}\sigma_{22} |\gamma_s|^2 & \sigma_{21}\gamma_s \bar\nu_s\\ 
       \sigma_{12} \bar \gamma_s \nu_s & \sigma_{11} |\nu_s|^2 \end{pmatrix}.
$$
By direct examination of the equation $\L_u(\sigma)=\lambda_u \,\sigma$ we obtain
\begin{equation}\label{eq_lambdauirrap}
\lambda_u=\sqrt{\ee(\exp\braket u A)\,}\sqrt{ \ee(\exp\braket uB)\,}.
\end{equation}
We immediately deduce
\[\lambda_u' = \braket u {\frac12(\ee(A)+\ee(B))}\qquad \lambda_u''-\lambda_u'{}^2 = \braket u {\frac12(\mathrm{var}\,A+\mathrm{var}\,B)\,u}.\]

In situation 2, we can use the extension discussed at the end of section~\ref{section_cltldp} with $P_{\R}=\ketbra{e_1}{e_1}$, and apply the formulas of Theorem \ref{theo_clt} with $\L$ replaced by~$\L_{\mathcal I_1(\cc e_1)}$. We see easily that the largest eigenvalue of $\L_u$ is 
\begin{equation}\label{eq_lambdausit2}
\lambda_u=\max(\sum_{s\in S}\e^{\braket us} |\alpha_s|^2, \sum_{s\in S}\e^{\braket us} |\beta_s|^2)
\end{equation}
and in a neighbourhood of zero, the first term is the largest, so that
\[\lambda_u'=\sum_{s\in S}\braket us\, |\alpha_s|^2 \quad \mbox{and}\quad \lambda_u'' =  \sum_{s\in S} {\braket u s}^2\, |\alpha_s|^2.\]

In situation 3, we again use the extension discussed at the end of section~\ref{section_cltldp} with ${\R_1}=\cc{e_1}$ and $\R_2=\cc e_2$.
The limit parameters for each corresponding restriction are computed in the previous point and correspond to those for the random variables $A$ and $B$. Since for any initial state $\rho$, a probability $\pp(X_0~=~i_0, \ldots, X_n=i_n)$ equals
\[\braket{e_1}{\rho(i_0)\,e_1}\, \prod_{k=1}^{n}|\alpha_{i_k-i_{k-1}}|^2+\braket{e_2}{\rho(i_0)\,e_2}\, \prod_{k=1}^{n}|\beta_{i_k-i_{k-1}}|^2\]
and we recover the parameters given in the statement above.
\fin
\begin{remark}
The irreducible periodic case described above can be understood in terms of a classical random walk, in a similar way to situation 3. Indeed, call a site $i$ in $V$ odd or even depending on the parity of its distance to the origin. Then exchanging the order of the basis vectors $e_1$ and $e_2$ at odd sites only is equivalent to considering a non-homogeneous OQRW with 
\[L_{i,i+s}=\begin{pmatrix}\nu_s&0\\ 0&\gamma_s\end{pmatrix}\quad \mbox{ if } i \mbox{ is even},\qquad
L_{i,i+s}=\begin{pmatrix}\gamma_s&0\\ 0&\nu_s\end{pmatrix}\quad \mbox{ if } i \mbox{ is odd}\]
(strictly speaking, such OQRWs do not enter into the framework of this article, but in the general case studied in \cite{CP1}).
Then, we define $(A_p)_{p\in \nn}$ and $(B_p)_{p\in \nn}$ to be two i.i.d. sequences with same law as $A$, $B$ respectively, and, if for example $X_0=0$ is even, we define a random variable $\pi$ to take the values $1$ and $2$ with probabilities $p=\braket{e_1}{\rho(i_0)\,e_1}$, $1-p$ respectively. Then, conditioned on $\pi=1$, the variable $X_p-X_0$ has the same law as $A_1+B_2+A_3+\ldots$ (where the sum stops at step~$p$). This explains the formulas given in Theorem~\ref{theo_hisc2} for situation~1, with $\L$ periodic, as well as the next proposition.
\end{remark}

\smallskip
For the case of irreducible, periodic $\L$ we also have a simpler explicit formula for the rate function of large deviations:
\begin{lemma}\label{lemma_ldpperiodicirr}
Assume an open quantum random walk with $V=\zz^d$ and $\h=\cc^2$ satisfies assumptions \emph{H1}, \emph{H2} and is irreducible periodic. Then, with the same notation as in Theorem \ref{theo_hisc2}, the position process $(X_p-X_0)_p$ satisfies a full large deviation principle, with rate function
\[c(u)=\frac12(\log\ee(\exp\braket uA)) + \log\ee(\exp\braket uB))).\]
\end{lemma}

\pre

This follows immediately from Theorem \ref{theo_ldp} and equation \eqref{eq_lambdauirrap} giving $\lambda_u$.
\fin
\begin{remark}\label{remark_LDPbreakdown}
In situation 2 of Lemma \ref{lemma_diag}, one sees that the largest eigenvalue $\lambda_u$ is given by \eqref{eq_lambdausit2}. For $u$ in a neighbourhood of zero, one has $\|\alpha_u\|>\|\beta_u\|$, but, if there exists $u$ such that $\|\alpha_u\|=\|\beta_u\|$, then $\lambda_u$ may not be differentiable and the large deviations principle may break down: see Example \ref{ex_ldpbreakdown}. A similar phenomenon can also appear in situation 3.
\end{remark}

%
%\begin{lemma}\label{lemma_irrper}
%Consider the irreducible, periodic case, as described above. Define $(A_p)_{p\in \nn}$ and $(B_p)_{p\in \nn}$ to be two processes of i.i.d. random variables with laws
%\[
%\pp(A_p=s)=|\nu_s|^2 \qquad \pp(B_p=s)=|\gamma_s|^2,
%\]
%by $m_A$, $m_B$ the expectations of $A$, $B$, by $C_A$, $C_B$ their covariance matrices and by $c_A$, $c_B$ their cumulant generating functions. The process $(X_p)_{p\in\nn}$ (or equivalently $(Q_p)_{p\in\nn}$) satisfies a central limit theorem with mean and covariance
%\[m=\frac12(m_A+m_B), \quad C =\frac12(C_A+C_B)\]
%and a large deviation property with a good rate function given as the Legendre transform of $c=\frac12(c_A+c_B)$. \end{lemma}
%
%\pre
%
%{ I don't find anymore $\eta=0$...}
%and one has
%\[\sum_{s\in S}\braket us^2\,\frac{|\alpha_s|^2+|\beta_s|^2}2=\tr\big(\L_u'(\rhoinv)\big)+2\,\tr\big(\L_u'(\eta_u)\big).\]
%Last, a straightforward computation shows that the eigenvalue with maximum modulus is $\lambda_u\|\gamma_u\|\,\|\nu_u\|$ so that 
%\[\log\lambda_u= \frac12\big(\log(\sum_{s\in S}\e^{\braket us}|\gamma_s|^2)+\log(\sum_{s\in S}\e^{\braket us}|\nu_s|^2)\big)=\frac12\big(c_A(u)+c_B(u)\big).\]
%Therefore, the formulas given in section \ref{section_cltldp} hold.
%\fin

\section{Examples}\label{section_examples}

%%%%%%%%%%%%%%%%%%%%
\begin{example}
We consider the case $d=1$, $\h={\mathbb C}^2$, $S=\{-1,+1\}$. In this case we can characterize irreducibility and period of the open quantum random walk $\M$ through the transition matrices $L_s, \, s=\pm 1$. In this example, we denote $L_-=L_{-1}$, $L_+=L_{+1}$. We state the next two propositions without proofs, as these are lengthy. The extension of these statements to finite homogeneous open quantum random walks, as well as the proofs, will be given in a future note.

\begin{prop}\label{reducibility-C2}
{\rm Irreducibility.}
Define 
\[W\overset{\mathrm{def}}{=} \{ \mbox{common eigenvectors of $L_+L_-$ and $L_-L_+$}\}.
\]
The homogeneous OQRW on ${\mathbb Z}$ is reducible if and only if one of the following facts holds
\begin{itemize}
\item $W$ contains an eigenvector of $L_-$ or $L_+$ %(this is automatically true when~$W~=~\mathbb C^2$), 
\item $W=\cc e_0\, \cup\, \cc e_1\setminus\{0\},$  for some linearly independent vectors $e_0$ and $e_1$ satisfying $L_-e_0,L_+e_0\in \cc e_1$ and $L_-e_1,L_+e_1\in \cc e_0$.
\end{itemize}
\end{prop}

\begin{prop} {\rm Period.} Suppose that the open quantum random walk $\M$ is irreducible. Its period can only be $2$ or $4$. It is $4$ if and only if there exists an orthonormal basis of ${\mathbb C}^2$ such that the representation of the transition matrices in that basis is
$$
 L_\varepsilon=\begin{pmatrix} a& 0\\ 0& b \end{pmatrix}, \qquad
 L_{-\varepsilon}=\begin{pmatrix}  0 & c \\  d & 0 \end{pmatrix}
 $$
for some $\varepsilon\in\{+,-\}$, where $a,b,c,d\in {\mathbb C}\setminus \{0\}$ are such that $|a|^2+|d|^2=|b|^2+|c|^2=1$.
 \end{prop}

\end{example}

\begin{example}\label{ex_stdexample}
We consider the standard example from \cite{APSS}, which is treated in section 5.3 of \cite{AGS}. This open quantum random walk is defined by $V=\zz$, $\h=\cc^2$, and transition operators given in the canonical basis $e_1$, $e_2$ of $\cc^2$ by
\[L_+=\frac1{\sqrt 3}\begin{pmatrix}1 & 1 \\ 0 & 1\end{pmatrix}\qquad L_-=\frac1{\sqrt 3}\begin{pmatrix}1 & 0 \\ -1 & 1\end{pmatrix}.\]
The only eigenvector of $L_+$ is $e_1$, the only eigenvector of $L_-$ is $e_2$, so that we are in situation 1 of Proposition \ref{prop_caracL} and $\L$ is irreducible. Again $L_+^2$ and $L_-^2$ have no eigenvector in common, so by Lemma \ref{lemma_Laperiodic}, we conclude that $\L$ is aperiodic (and therefore regular, by Lemma \ref{lemma_irrapereg}). We observe that $\rhoinv=\frac12\id$ is the invariant state of $\L$. We compute the quantities $m$ and $C\in\rr_+$ from Theorem \ref{theo_clt}:
\[m= \tr(L_+ L_+^*)-\tr(L_- L_-^*)=0.\] 
To compute $C$ we need to find the solution $\eta$ of 
\[(\id -\L)(\eta)= \frac16 \begin{pmatrix} \hphantom{,}1 & \hphantom{+}2 \\ \hphantom{,}2 & -1\end{pmatrix}\]
satisfying $\tr\,\eta=0$.
We find $\eta=\frac1{12}\begin{pmatrix} 5 & \hphantom{-}2 \\ 2 & -5\end{pmatrix}$, and we have
\begin{equation*}
C= \tr\big(L_+ \rhoinv L_+^*+ L_-\rhoinv L_-^*) + 2\,\tr\big(L_+ \eta L_+^*-L_-\eta L_-^*\big)=\frac89.
\end{equation*}
By Theorem \ref{theo_clt}, we have the convergence in law
\[\frac{X_p-X_0}{\sqrt p}\underset{p\to\infty}{\longrightarrow} \mathcal N(0, \frac89).\]
To verify the validity of this statement, in Figure \ref{figure_CLT1} we display the empirical cumulative distribution function of a 1000-sample of $\frac{X_p}{\sqrt{8p/9}}$ conditioned on~$X_0=~0$ for $p=10,100,1000$, and compare it to the cumulative distribution function of a standard normal variable.
\begin{remark}
In this and the following simulations, the initial state is assumed to be of the form $\rho(0)\otimes\ketbra 00$ and $\rho(0)$ is chosen randomly as $\frac{XX'}{\tr(XX')}$ where $X$ has independent entries with uniform law on $[0,1]$. 
\end{remark}
 
\begin{figure*}[h]
%\hspace{-2cm}
\begin{center}
\includegraphics[width=1\textwidth]{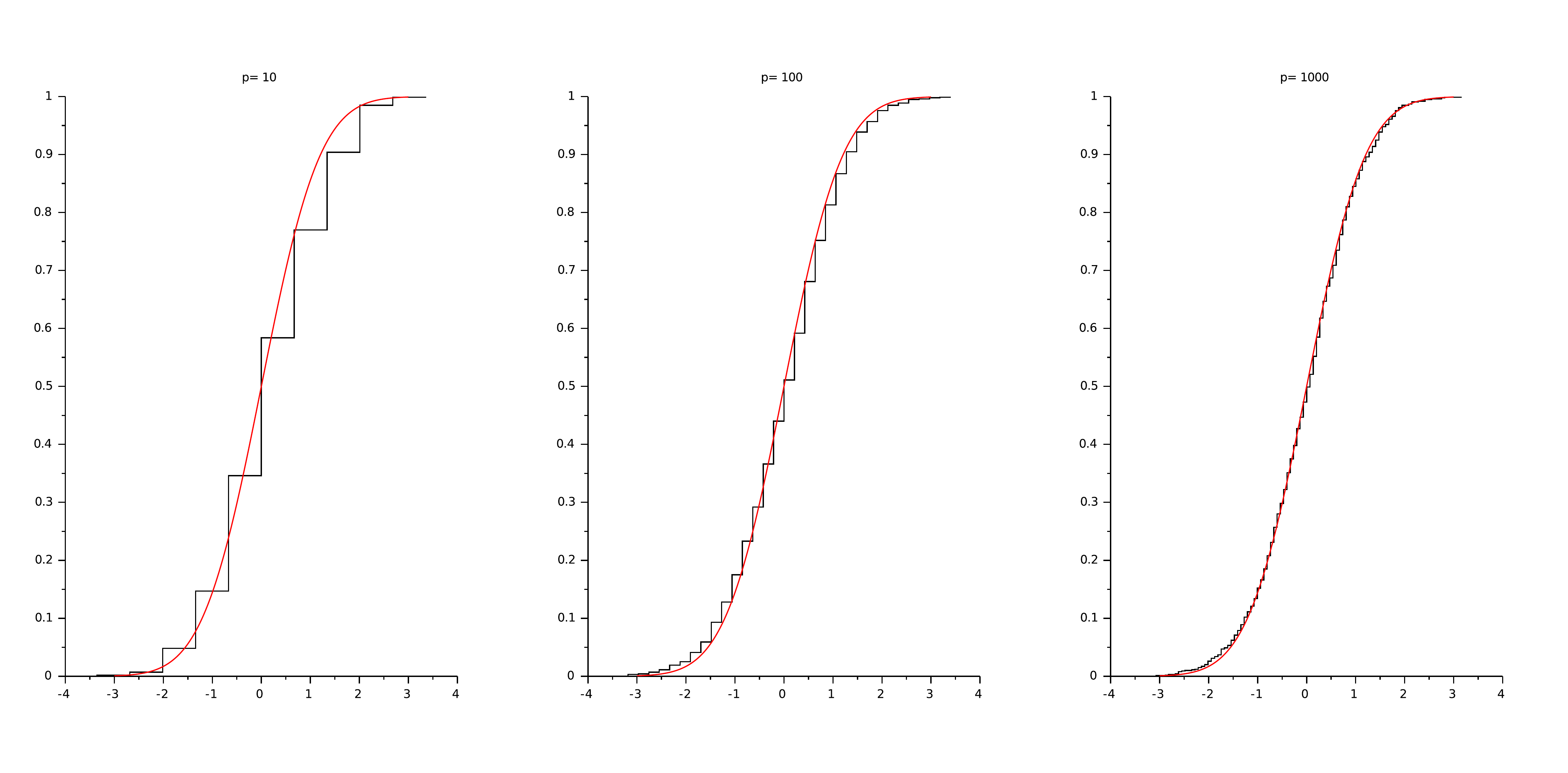}
\caption{C.D.F. of a 1000-sample of $(\frac{X_p}{\sqrt{8p/9}})_p$ with $X_0=0$, in Example \ref{ex_stdexample}.}
\label{figure_CLT1}
\end{center}
\end{figure*}

By Theorem \ref{theo_ldp}, the process $(\frac{X_p-X_0}p)_p$ satisfies a large deviation property with good rate function equal to the Legendre transform $I$ of $u\mapsto \log\lambda_u$, where $\lambda_u$ is the largest eigenvalue of $\L_u$. This map $\L_u$, written in the canonical basis of the set of two by two matrices, has basis 
\[\frac13 \begin{pmatrix}
\e^{u}+\e^{-u} & \e^{u} & \e^{u} & \e^{u} \\
-\e^{-u} & \e^{u}+\e^{-u} & 0 & \e^{u}\\
-\e^{-u} & 0 & \e^{u}+\e^{-u} & \e^{u}\\
\e^{-u} & -\e^{-u} & -\e^{-u} & \e^{u}+\e^{-u}
\end{pmatrix} \]
and by a tedious computation, one shows that $\lambda_u$ equals
\[\frac13 \big(\e^{u}+\e^{-u}+(\e^{u}+\e^{-u}+\sqrt{\e^{2u}+\e^{-2u}+3})^{1/3} - (\e^{u}+\e^{-u}+\sqrt{\e^{2u}+\e^{-2u}+3})^{-1/3}\big).\]
As expected from Lemma \ref{lemma_analyticitylambdau}, this is a smooth and strictly convex function. Numerical computations prove that the rate function $I$ has the form displayed in Figure \ref{figure_ratefunction}.
\begin{figure*}[h!]
\begin{center}
\includegraphics[width=0.6\textwidth]{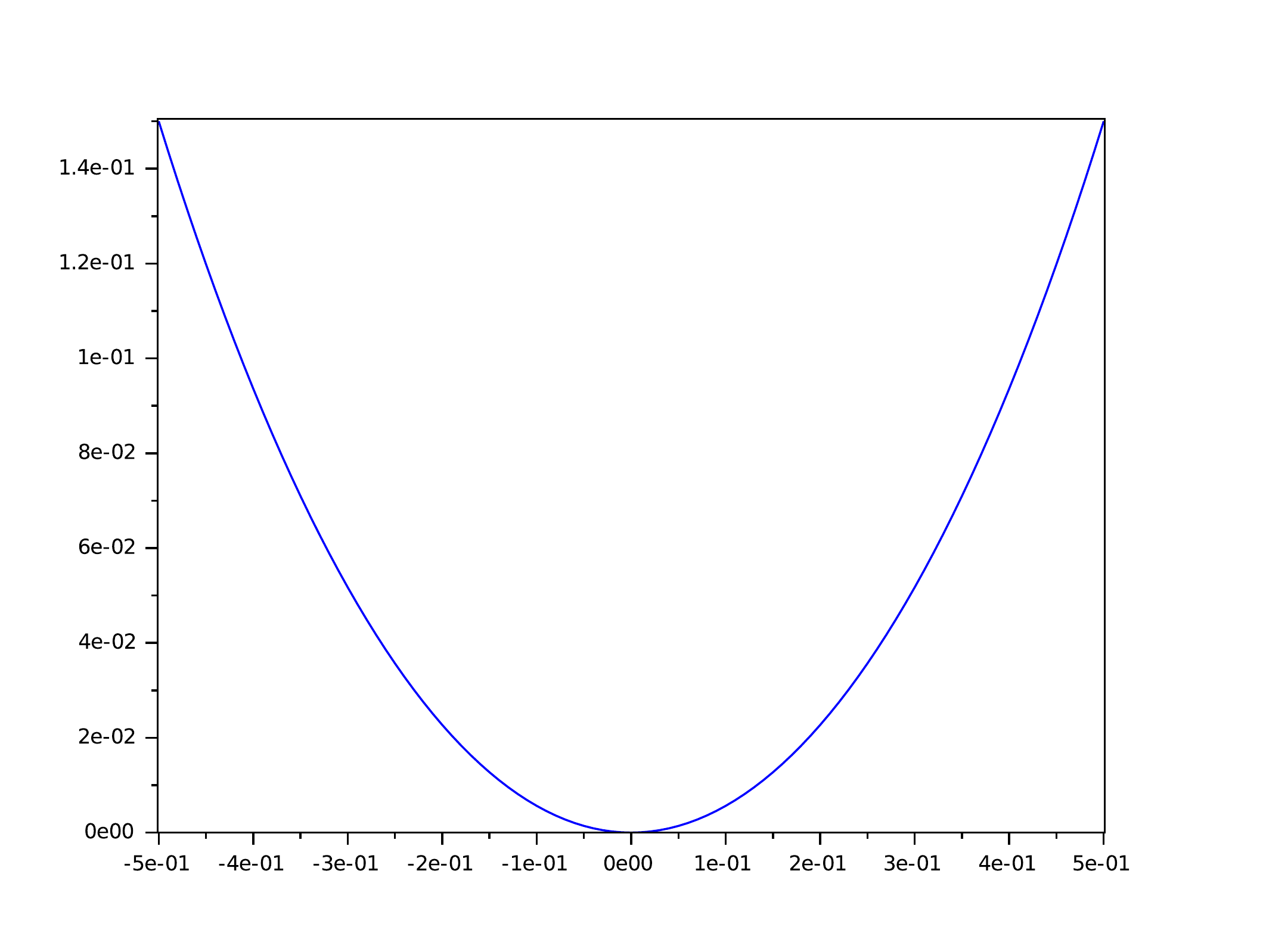}
\caption{Rate function for $(\frac{X_p-X_0}p)_p$ in Example \ref{ex_stdexample} }
\label{figure_ratefunction}
\end{center}
\end{figure*}
\end{example}

\begin{example}\label{ex_Lnonirr}
We consider the open quantum random walk defined by $V=\zz$, $\h=\cc^2$, and transition operators given in the canonical basis $e_1$, $e_2$ of $\cc^2$ by
\[L_+=\begin{pmatrix}0 & \sqrt 3 /2 \\1/\sqrt2  & 0\end{pmatrix}\qquad L_-=\begin{pmatrix}0 & 1/2 \\ 1/\sqrt 2 & 0\end{pmatrix}.\]
From Lemma \ref{lemma_Laperiodic}, the map $\L$ is irreducible and $2$-periodic. Then according to Theorem \ref{theo_hisc2}, defining $A$ and $B$ to be random variables with values in $S$ satisfying
\[\pp(A=+1)=1/2,\qquad \pp(A=-1)=1/2,\]
\[\pp(B=+1)=3/4,\qquad \pp(B=-1)=1/4,\]
with mean, variance, and cumulant generating function
\[m_A=0,\quad C_A=1,\quad c_A(u)=\log(\e^{u}+\e^{-u})-\log 2\]
\[m_B=1/2,\quad C_B=3/4,\quad c_B(u)=\log(3\e^{u}+\e^{-u})-2\log 2\]
then with the notation of Theorem \ref{theo_hisc2} 
\[m=(m_A+m_B)/2 = 1/4\qquad C=(C_A+C_B)/2=7/8\]
and one has the convergence in law
\[\frac{X_p- p/4}{\sqrt p}\underset{p\to\infty}{\longrightarrow} \mathcal N(0,\frac78).\]
Figure \ref{figure_CLT3} below displays the cumulative distribution function of a 1000-sample of $\frac{X_p-p/4}{\sqrt{7p/8}}$, conditioned on $X_0=0$, and the cumulative distribution function of a standard normal variable for $p=10,100,1000$.
\begin{figure*}[h!]
%\hspace{-2cm}
\begin{center}
\includegraphics[width=1.0\textwidth]{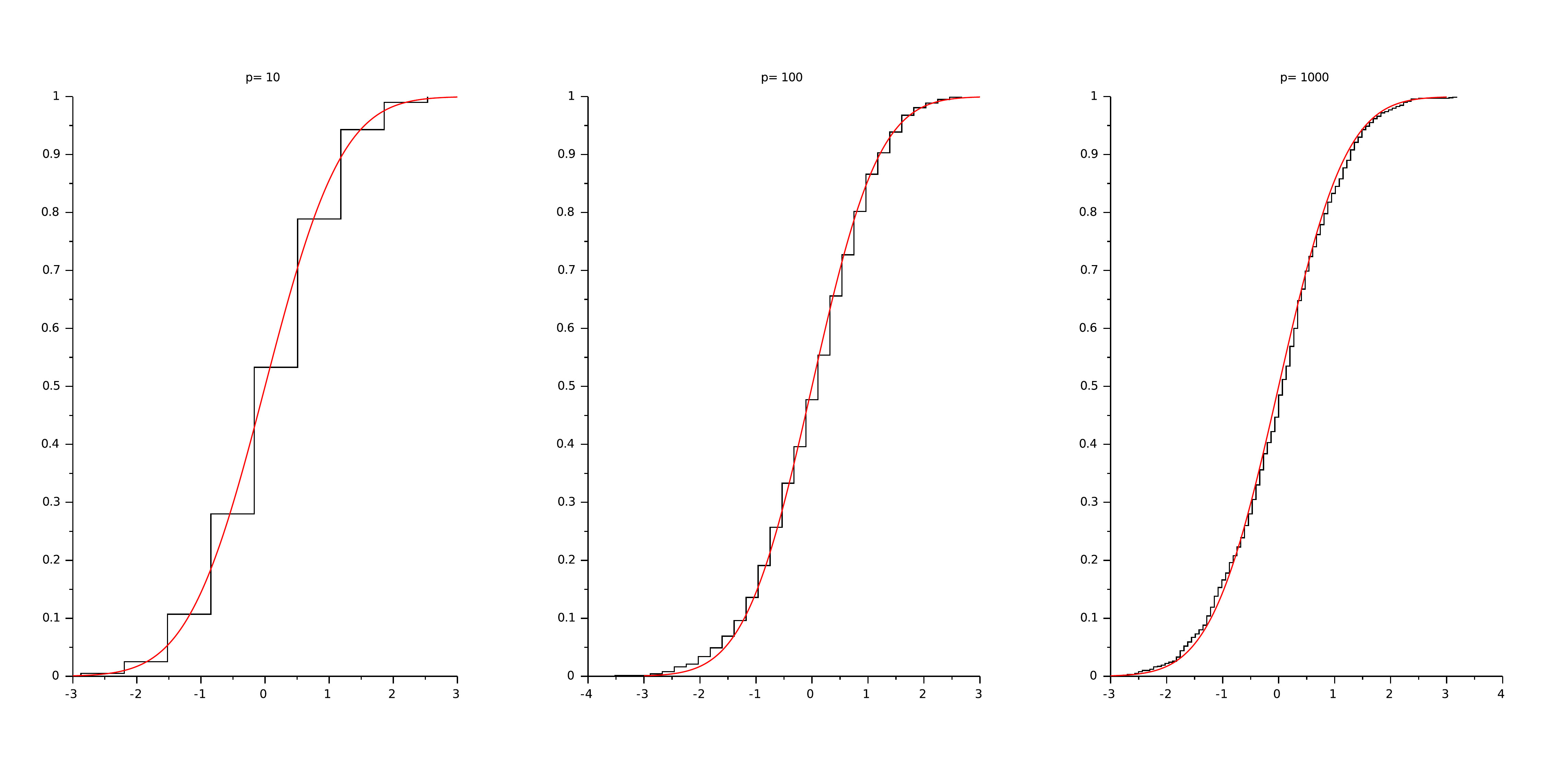}
\caption{C.D.F. of a 1000-sample of $\frac{X_p-p/4}{\sqrt{7p/8}}$ with $X_0=0$, in Example \ref{ex_Lnonirr}.}
\label{figure_CLT3}
\end{center}
\end{figure*}

In addition, the process $(\frac{X_p-X_0}p)_{p\in\nn}$ satisfies a large deviation property with a good rate function $I$ obtained as the Legendre transform of
\[c(u)=\frac12(c_A(u)+c_B(u))= \frac12\big(\log(\e^{u}+\e^{-u})+\log(3\e^{u}+\e^{-u})\big)-\frac32\log 2.\]
Explicitly, one finds that $I(t)=+\infty$ for $t\not\in]-1,1[$ and for $t\in]-1,+1[$:
\[I(t)= t\, u_t + \frac32\log 2 - \frac12\big(\log(\e^{u_t}+\e^{-u_t})+\log(3\e^{u_t}+\e^{-u_t})\big)\]
where $u_t=\frac12\, \log\,\frac{2t+\sqrt{t^2+3}}{3(1-t)}$. This rate function has the profile displayed in Figure \ref{figure_Lnonirr}.
\begin{figure*}[h!]
\begin{center}
\includegraphics[width=0.6\textwidth]{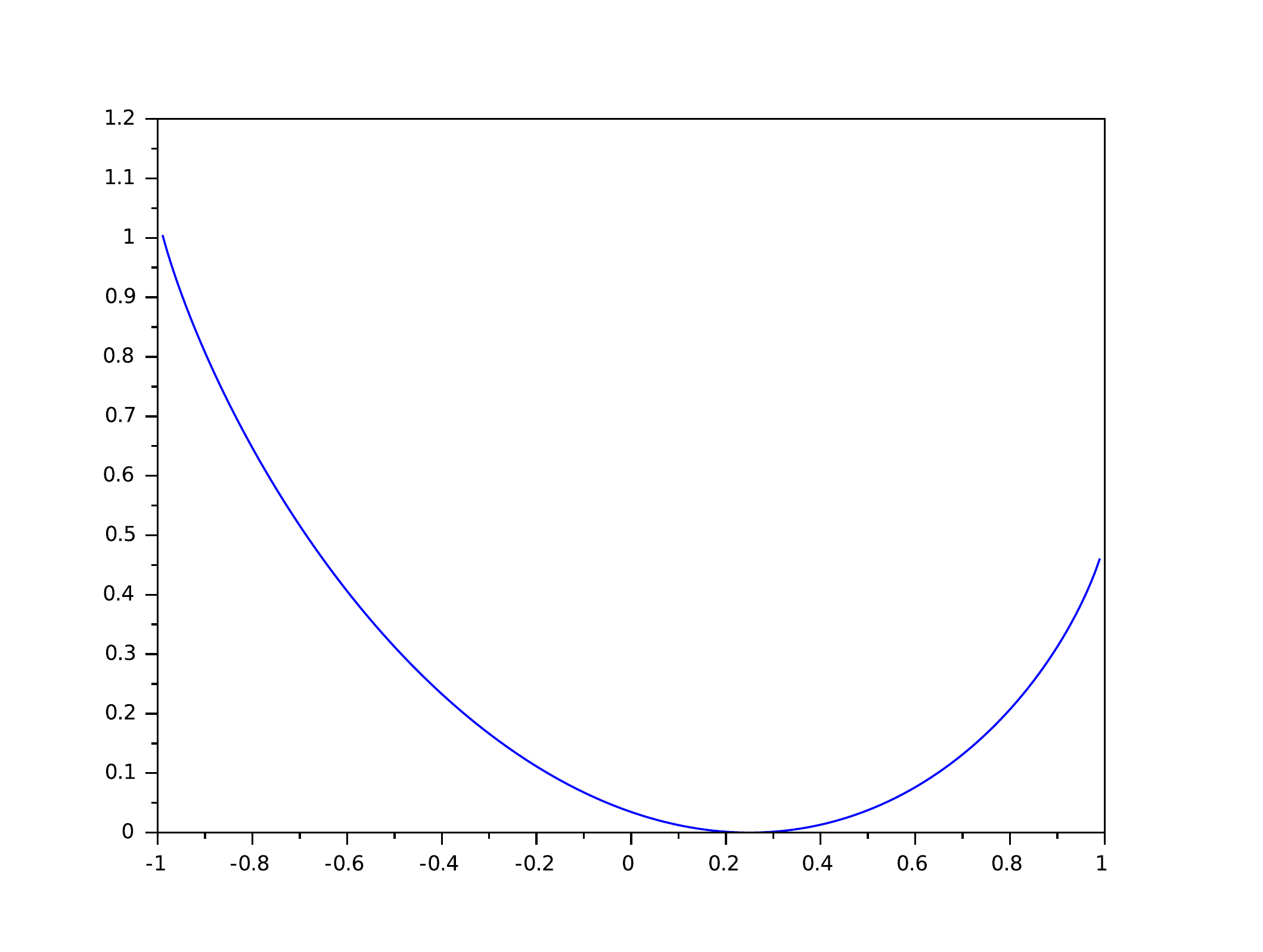}
\caption{Rate function for $(\frac{X_p-X_0}p)_p$ in Example \ref{ex_Lnonirr}}
\label{figure_Lnonirr}
\end{center}
\end{figure*}

\end{example}

\begin{example}\label{ex_ldpbreakdown}
Consider the open quantum random walk defined by $V=\zz$, $\h=\cc^2$, and transition operators given in the canonical basis $e_1$, $e_2$ of $\cc^2$ by
\[L_+=\begin{pmatrix}\frac1{\sqrt 2} & \frac1{2\sqrt 2} \\ 0 & \frac{\sqrt 3}2\end{pmatrix}\qquad L_-=\begin{pmatrix}\frac1{\sqrt 2} & -\frac1{2\sqrt 2} \\ 0 & 0 \end{pmatrix}.\]
First observe that the map $\L$ is not irreducible in this case, as we are in situation~2 of Proposition \ref{prop_caracL}. A straightforward computation shows that the largest eigenvalue of $\L_u$ is 
\[\lambda_u=\sup(\frac{\e^u+\e^{-u}}2, \frac {3\,\e^u}{4}).\]
For $u$ close to zero $\lambda_u$ is $\frac{\e^u+\e^{-u}}2$ so that $\lambda_u'=0$ and $\lambda_u''=1$ for $u=0$. We must therefore have
\[\frac{X_p-X_0}{\sqrt p} \underset{p\to\infty}{\longrightarrow} \mathcal N(0,1).\]
Figure \ref{figure_CLT4} below displays the cumulative distribution function of a 1000-sample of $\frac{X_p}{\sqrt{p}}$, conditioned on $X_0=0$, and the cumulative distribution function of a standard normal variable for $p=10,100,1000$.
\begin{figure*}[h]
%\hspace{-2cm}
\begin{center}
\includegraphics[width=1.0\textwidth]{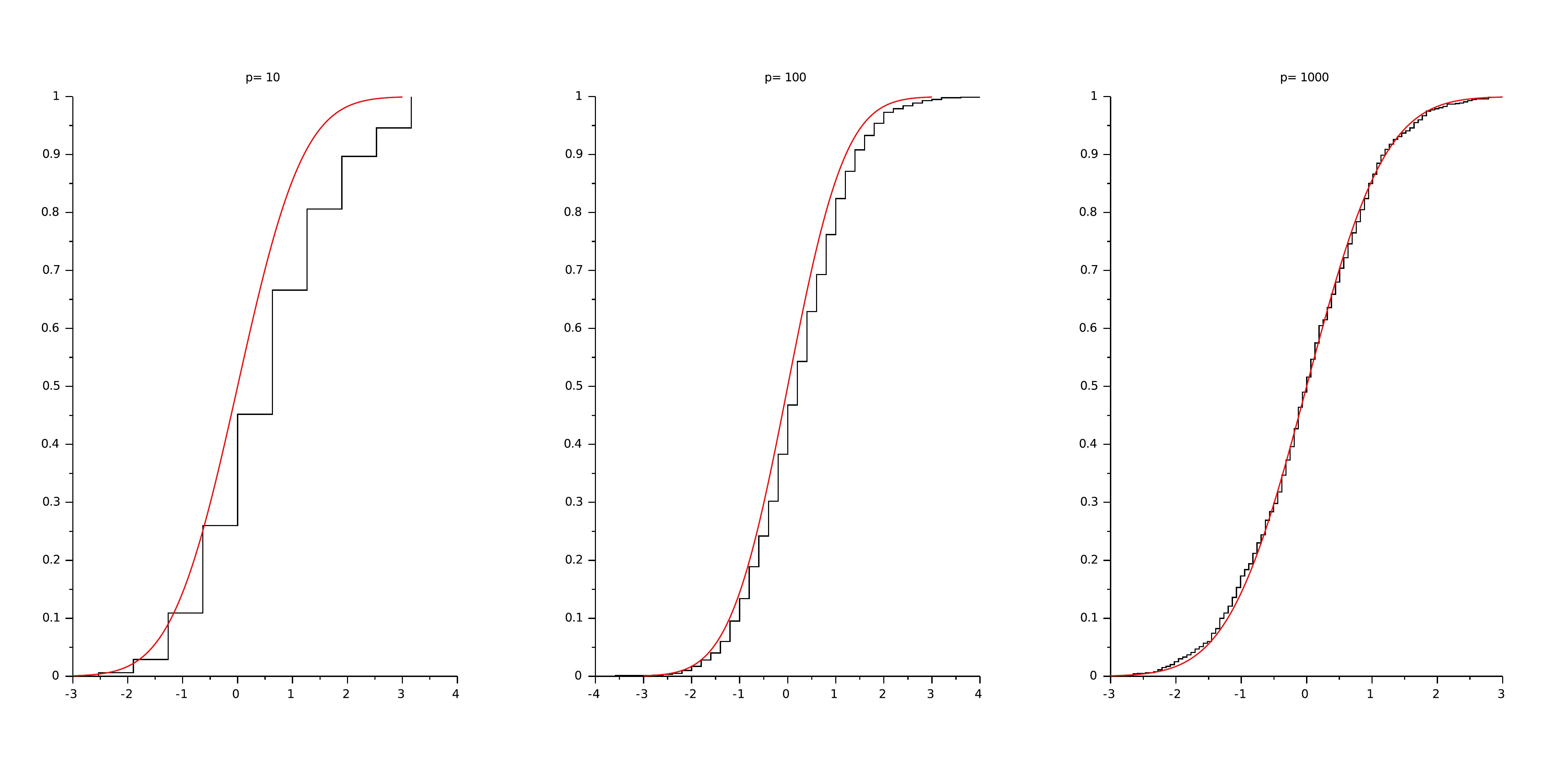}
\caption{C.D.F. of a 1000-sample of $(\frac{X_p-X_0}{\sqrt{p}})_p$ with $X_0=0$, in Example \ref{ex_ldpbreakdown}.}
\label{figure_CLT4}
\end{center}
\end{figure*}

Due to the generalizations discussed at the end of Section \ref{section_cltldp}, we have $\R=\cc e_1$, $\D=\cc e_2$ and the central limit theorem holds: the behavior of the process~$(X_p)_p$, associated with $\L$, is the same as the one of the process $(\widetilde X_p)_p$ associated with the restriction $\L_{|\R}$.  

As we commented previously, giving a large deviations result in this case is harder and we cannot use the general results we proved. A G\"artner-Ellis theorem could be applied by direct computation of the moment generating functions. In general, however, the rate function for the process $(X_p)_p$ will not coincide with the one for $(\widetilde X_p)_p$, since it will essentially depend on how much time the evolution spends in $\D$.

More precisely, for the transition matrices introduced above and taking the initial state $\rho=\ketbra {e_2}{e_2}\otimes \ketbra 0 0$, we have, by relation \eqref{eq_probatraj},
$$
P(X_n = n) = \tr(\ketbra {L_+^ne_2}{L_+^ne_2})
= \left(\frac{3}{4}\right)^n + \left(\frac{1}{8}\right)2^{1-n} \frac{\sqrt {3^n}-\sqrt {2^n}}{\sqrt 3- \sqrt 2} 
$$
and consequently 
$$\lim_n \frac 1n \log \ee[e^{uX_n}] \ge \log\big(\frac 34 e^u\big)\qquad \mbox{ for all } u,$$
while $\lim_n \frac 1n \log \ee[e^{u\widetilde X_n}] = \log\big(\frac {e^u+e^{-u}}{2}\big)$, which for $u>\log 2$ is smaller than the bound 
$\log\big(\frac 34 e^u\big)$.

This clarifies the fact that the large deviations will not depend only on $\L_{|\R}$. Moreover, a second problem arises in this example, which is the lack of regularity of $\lambda_u$. Indeed, $\lambda_u$ is the supremum of two quantities which coincide for $u_0=\frac12\log2$, and $\log\lambda_u$ is not differentiable at $u_0$: the left derivative is equal to $\frac{\sqrt2}4$ and the right derivative to $\frac{3\sqrt2}4$. The Legendre transform of $\log\lambda_u$ is displayed in Figure \ref{figure_ldpbreakdown}, and we observe that it is not strictly convex.
\begin{figure*}[h]
\begin{center}
\includegraphics[width=0.6\textwidth]{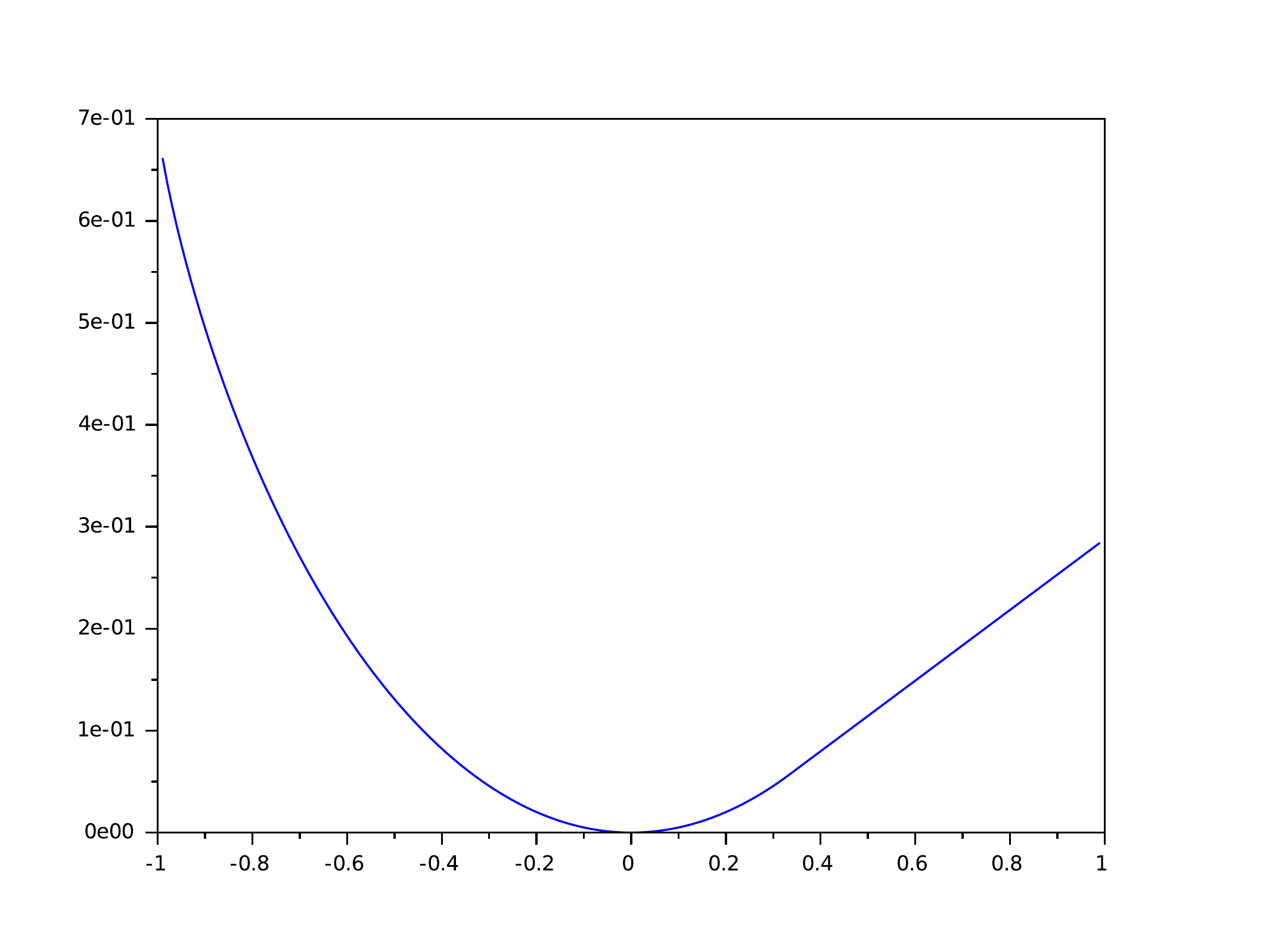}
\caption{Rate function for $(\frac{X_p-X_0}p)_p$ in Example \ref{ex_ldpbreakdown}}
\label{figure_ldpbreakdown}
\end{center}
\end{figure*}
\end{example}

\bibliography{biblio}

\end{document}